\input amstex
%

\def\next{AMS-SEKR}\ifx\styname\next \endinput\fi
\catcode`\@=11
\def\styname{AMS-SEKR}
\def\styversion{2.0}
{\W@{}\W@{\styname.STY - Version \styversion}\W@{}}
\hyphenation{acad-e-my acad-e-mies af-ter-thought anom-aly anom-alies
an-ti-deriv-a-tive an-tin-o-my an-tin-o-mies apoth-e-o-ses apoth-e-o-sis
ap-pen-dix ar-che-typ-al as-sign-a-ble as-sist-ant-ship as-ymp-tot-ic
asyn-chro-nous at-trib-uted at-trib-ut-able bank-rupt bank-rupt-cy
bi-dif-fer-en-tial blue-print busier busiest cat-a-stroph-ic
cat-a-stroph-i-cally con-gress cross-hatched data-base de-fin-i-tive
de-riv-a-tive dis-trib-ute dri-ver dri-vers eco-nom-ics econ-o-mist
elit-ist equi-vari-ant ex-quis-ite ex-tra-or-di-nary flow-chart
for-mi-da-ble forth-right friv-o-lous ge-o-des-ic ge-o-det-ic geo-met-ric
griev-ance griev-ous griev-ous-ly hexa-dec-i-mal ho-lo-no-my ho-mo-thetic
ideals idio-syn-crasy in-fin-ite-ly in-fin-i-tes-i-mal ir-rev-o-ca-ble
key-stroke lam-en-ta-ble light-weight mal-a-prop-ism man-u-script
mar-gin-al meta-bol-ic me-tab-o-lism meta-lan-guage me-trop-o-lis
met-ro-pol-i-tan mi-nut-est mol-e-cule mono-chrome mono-pole mo-nop-oly
mono-spline mo-not-o-nous mul-ti-fac-eted mul-ti-plic-able non-euclid-ean
non-iso-mor-phic non-smooth par-a-digm par-a-bol-ic pa-rab-o-loid
pa-ram-e-trize para-mount pen-ta-gon phe-nom-e-non post-script pre-am-ble
pro-ce-dur-al pro-hib-i-tive pro-hib-i-tive-ly pseu-do-dif-fer-en-tial
pseu-do-fi-nite pseu-do-nym qua-drat-ics quad-ra-ture qua-si-smooth
qua-si-sta-tion-ary qua-si-tri-an-gu-lar quin-tes-sence quin-tes-sen-tial
re-arrange-ment rec-tan-gle ret-ri-bu-tion retro-fit retro-fit-ted
right-eous right-eous-ness ro-bot ro-bot-ics sched-ul-ing se-mes-ter
semi-def-i-nite semi-ho-mo-thet-ic set-up se-vere-ly side-step sov-er-eign
spe-cious spher-oid spher-oid-al star-tling star-tling-ly
sta-tis-tics sto-chas-tic straight-est strange-ness strat-a-gem strong-hold
sum-ma-ble symp-to-matic syn-chro-nous topo-graph-i-cal tra-vers-a-ble
tra-ver-sal tra-ver-sals treach-ery turn-around un-at-tached un-err-ing-ly
white-space wide-spread wing-spread wretch-ed wretch-ed-ly Brown-ian
Eng-lish Euler-ian Feb-ru-ary Gauss-ian Grothen-dieck Hamil-ton-ian
Her-mit-ian Jan-u-ary Japan-ese Kor-te-weg Le-gendre Lip-schitz
Lip-schitz-ian Mar-kov-ian Noe-ther-ian No-vem-ber Rie-mann-ian
Schwarz-schild Sep-tem-ber
form per-iods Uni-ver-si-ty cri-ti-sism for-ma-lism}
\Invalid@\nofrills
\Invalid@\usualspace
\newif\ifnofrills@
\def\nofrills@#1#2{\relaxnext@
  \DN@{\ifx\next\nofrills
    \nofrills@true\let#2\relax\DN@\nofrills{\nextii@}%
  \else
    \nofrills@false\def#2{#1}\let\next@\nextii@\fi
\next@}}
\def\usualspace@#1{\ifnofrills@\def\usualspace{#1}\fi}
\def\addto#1#2{\csname \expandafter\eat@\string#1@\endcsname
  \expandafter{\the\csname \expandafter\eat@\string#1@\endcsname#2}}
\newdimen\bigsize@
\def\big@#1#2{{\hbox{$\left#2\vcenter to#1\bigsize@{}%
  \right.\nulldelimiterspace\z@\m@th$}}}
\def\big{\big@\@ne}
\def\Big{\big@{1.5}}
\def\bigg{\big@\tw@}
\def\Bigg{\big@{2.5}}
\def\raggedcenter@{\leftskip\z@ plus.4\hsize \rightskip\leftskip
 \parfillskip\z@ \parindent\z@ \spaceskip.3333em \xspaceskip.5em
 \pretolerance9999\tolerance9999 \exhyphenpenalty\@M
 \hyphenpenalty\@M \let\\\linebreak}
\def\upperspecialchars{\def\ss{SS}\let\i=I\let\j=J\let\ae\AE\let\oe\OE
  \let\o\O\let\aa\AA\let\l\L}
\def\uppercasetext@#1{%
  {\spaceskip1.2\fontdimen2\the\font plus1.2\fontdimen3\the\font
   \upperspecialchars\uctext@#1$\m@th\aftergroup\eat@$}}
\def\uctext@#1$#2${\endash@#1-\endash@$#2$\uctext@}
\def\endash@#1-#2\endash@{\uppercase{#1}\if\notempty{#2}--\endash@#2\endash@\fi}
\def\runaway@#1{\DN@{#1}\ifx\envir@\next@
  \Err@{You seem to have a missing or misspelled \string\end#1 ...}%
  \let\envir@\empty\fi}
\newif\iftemp@
\def\notempty#1{TT\fi\def\test@{#1}\ifx\test@\empty\temp@false
  \else\temp@true\fi \iftemp@}
\font@\tensmc=cmcsc10
\font@\sevenex=cmex7
\font@\sevenit=cmti7
\font@\eightrm=cmr8 
\font@\sixrm=cmr6 
\font@\eighti=cmmi8     \skewchar\eighti='177 
\font@\sixi=cmmi6       \skewchar\sixi='177   
\font@\eightsy=cmsy8    \skewchar\eightsy='60 
\font@\sixsy=cmsy6      \skewchar\sixsy='60   
\font@\eightex=cmex8
\font@\eightbf=cmbx8 
\font@\sixbf=cmbx6   
\font@\eightit=cmti8 
\font@\eightsl=cmsl8 
\font@\eightsmc=cmcsc8
\font@\eighttt=cmtt8 


\loadmsam
\loadmsbm
\loadeufm
\UseAMSsymbols
\newtoks\tenpoint@
\def\tenpoint{\normalbaselineskip12\p@
 \abovedisplayskip12\p@ plus3\p@ minus9\p@
 \belowdisplayskip\abovedisplayskip
 \abovedisplayshortskip\z@ plus3\p@
 \belowdisplayshortskip7\p@ plus3\p@ minus4\p@
 \textonlyfont@\rm\tenrm \textonlyfont@\it\tenit
 \textonlyfont@\sl\tensl \textonlyfont@\bf\tenbf
 \textonlyfont@\smc\tensmc \textonlyfont@\tt\tentt
 \textonlyfont@\bsmc\tenbsmc
 \ifsyntax@ \def\big##1{{\hbox{$\left##1\right.$}}}%
  \let\Big\big \let\bigg\big \let\Bigg\big
 \else
  \textfont\z@=\tenrm  \scriptfont\z@=\sevenrm  \scriptscriptfont\z@=\fiverm
  \textfont\@ne=\teni  \scriptfont\@ne=\seveni  \scriptscriptfont\@ne=\fivei
  \textfont\tw@=\tensy \scriptfont\tw@=\sevensy \scriptscriptfont\tw@=\fivesy
  \textfont\thr@@=\tenex \scriptfont\thr@@=\sevenex
        \scriptscriptfont\thr@@=\sevenex
  \textfont\itfam=\tenit \scriptfont\itfam=\sevenit
        \scriptscriptfont\itfam=\sevenit
  \textfont\bffam=\tenbf \scriptfont\bffam=\sevenbf
        \scriptscriptfont\bffam=\fivebf
  \setbox\strutbox\hbox{\vrule height8.5\p@ depth3.5\p@ width\z@}%
  \setbox\strutbox@\hbox{\lower.5\normallineskiplimit\vbox{%
        \kern-\normallineskiplimit\copy\strutbox}}%
 \setbox\z@\vbox{\hbox{$($}\kern\z@}\bigsize@=1.2\ht\z@
 \fi
 \normalbaselines\rm\ex@.2326ex\jot3\ex@\the\tenpoint@}
\newtoks\eightpoint@
\def\eightpoint{\normalbaselineskip10\p@
 \abovedisplayskip10\p@ plus2.4\p@ minus7.2\p@
 \belowdisplayskip\abovedisplayskip
 \abovedisplayshortskip\z@ plus2.4\p@
 \belowdisplayshortskip5.6\p@ plus2.4\p@ minus3.2\p@
 \textonlyfont@\rm\eightrm \textonlyfont@\it\eightit
 \textonlyfont@\sl\eightsl \textonlyfont@\bf\eightbf
 \textonlyfont@\smc\eightsmc \textonlyfont@\tt\eighttt
 \textonlyfont@\bsmc\eightbsmc
 \ifsyntax@\def\big##1{{\hbox{$\left##1\right.$}}}%
  \let\Big\big \let\bigg\big \let\Bigg\big
 \else
  \textfont\z@=\eightrm \scriptfont\z@=\sixrm \scriptscriptfont\z@=\fiverm
  \textfont\@ne=\eighti \scriptfont\@ne=\sixi \scriptscriptfont\@ne=\fivei
  \textfont\tw@=\eightsy \scriptfont\tw@=\sixsy \scriptscriptfont\tw@=\fivesy
  \textfont\thr@@=\eightex \scriptfont\thr@@=\sevenex
   \scriptscriptfont\thr@@=\sevenex
  \textfont\itfam=\eightit \scriptfont\itfam=\sevenit
   \scriptscriptfont\itfam=\sevenit
  \textfont\bffam=\eightbf \scriptfont\bffam=\sixbf
   \scriptscriptfont\bffam=\fivebf
 \setbox\strutbox\hbox{\vrule height7\p@ depth3\p@ width\z@}%
 \setbox\strutbox@\hbox{\raise.5\normallineskiplimit\vbox{%
   \kern-\normallineskiplimit\copy\strutbox}}%
 \setbox\z@\vbox{\hbox{$($}\kern\z@}\bigsize@=1.2\ht\z@
 \fi
 \normalbaselines\eightrm\ex@.2326ex\jot3\ex@\the\eightpoint@}

\font@\twelverm=cmr10 scaled\magstep1
\font@\twelveit=cmti10 scaled\magstep1
\font@\twelvesl=cmsl10 scaled\magstep1
\font@\twelvesmc=cmcsc10 scaled\magstep1
\font@\twelvett=cmtt10 scaled\magstep1
\font@\twelvebf=cmbx10 scaled\magstep1
\font@\twelvei=cmmi10 scaled\magstep1
\font@\twelvesy=cmsy10 scaled\magstep1
\font@\twelveex=cmex10 scaled\magstep1
\font@\twelvemsa=msam10 scaled\magstep1
\font@\twelveeufm=eufm10 scaled\magstep1
\font@\twelvemsb=msbm10 scaled\magstep1
\newtoks\twelvepoint@
\def\twelvepoint{\normalbaselineskip15\p@
 \abovedisplayskip15\p@ plus3.6\p@ minus10.8\p@
 \belowdisplayskip\abovedisplayskip
 \abovedisplayshortskip\z@ plus3.6\p@
 \belowdisplayshortskip8.4\p@ plus3.6\p@ minus4.8\p@
 \textonlyfont@\rm\twelverm \textonlyfont@\it\twelveit
 \textonlyfont@\sl\twelvesl \textonlyfont@\bf\twelvebf
 \textonlyfont@\smc\twelvesmc \textonlyfont@\tt\twelvett
 \textonlyfont@\bsmc\twelvebsmc
 \ifsyntax@ \def\big##1{{\hbox{$\left##1\right.$}}}%
  \let\Big\big \let\bigg\big \let\Bigg\big
 \else
  \textfont\z@=\twelverm  \scriptfont\z@=\tenrm  \scriptscriptfont\z@=\sevenrm
  \textfont\@ne=\twelvei  \scriptfont\@ne=\teni  \scriptscriptfont\@ne=\seveni
  \textfont\tw@=\twelvesy \scriptfont\tw@=\tensy \scriptscriptfont\tw@=\sevensy
  \textfont\thr@@=\twelveex \scriptfont\thr@@=\tenex
        \scriptscriptfont\thr@@=\tenex
  \textfont\itfam=\twelveit \scriptfont\itfam=\tenit
        \scriptscriptfont\itfam=\tenit
  \textfont\bffam=\twelvebf \scriptfont\bffam=\tenbf
        \scriptscriptfont\bffam=\sevenbf
  \setbox\strutbox\hbox{\vrule height10.2\p@ depth4.2\p@ width\z@}%
  \setbox\strutbox@\hbox{\lower.6\normallineskiplimit\vbox{%
        \kern-\normallineskiplimit\copy\strutbox}}%
 \setbox\z@\vbox{\hbox{$($}\kern\z@}\bigsize@=1.4\ht\z@
 \fi
 \normalbaselines\rm\ex@.2326ex\jot3.6\ex@\the\twelvepoint@}

\def\headfonts{\twelvepoint\bf}

\font@\fourteenrm=cmr10 scaled\magstep2
\font@\fourteenit=cmti10 scaled\magstep2
\font@\fourteensl=cmsl10 scaled\magstep2
\font@\fourteensmc=cmcsc10 scaled\magstep2
\font@\fourteentt=cmtt10 scaled\magstep2
\font@\fourteenbf=cmbx10 scaled\magstep2
\font@\fourteeni=cmmi10 scaled\magstep2
\font@\fourteensy=cmsy10 scaled\magstep2
\font@\fourteenex=cmex10 scaled\magstep2
\font@\fourteenmsa=msam10 scaled\magstep2
\font@\fourteeneufm=eufm10 scaled\magstep2
\font@\fourteenmsb=msbm10 scaled\magstep2
\newtoks\fourteenpoint@
\def\fourteenpoint{\normalbaselineskip15\p@
 \abovedisplayskip18\p@ plus4.3\p@ minus12.9\p@
 \belowdisplayskip\abovedisplayskip
 \abovedisplayshortskip\z@ plus4.3\p@
 \belowdisplayshortskip10.1\p@ plus4.3\p@ minus5.8\p@
 \textonlyfont@\rm\fourteenrm \textonlyfont@\it\fourteenit
 \textonlyfont@\sl\fourteensl \textonlyfont@\bf\fourteenbf
 \textonlyfont@\smc\fourteensmc \textonlyfont@\tt\fourteentt
 \textonlyfont@\bsmc\fourteenbsmc
 \ifsyntax@ \def\big##1{{\hbox{$\left##1\right.$}}}%
  \let\Big\big \let\bigg\big \let\Bigg\big
 \else
  \textfont\z@=\fourteenrm  \scriptfont\z@=\twelverm  \scriptscriptfont\z@=\tenrm
  \textfont\@ne=\fourteeni  \scriptfont\@ne=\twelvei  \scriptscriptfont\@ne=\teni
  \textfont\tw@=\fourteensy \scriptfont\tw@=\twelvesy \scriptscriptfont\tw@=\tensy
  \textfont\thr@@=\fourteenex \scriptfont\thr@@=\twelveex
        \scriptscriptfont\thr@@=\twelveex
  \textfont\itfam=\fourteenit \scriptfont\itfam=\twelveit
        \scriptscriptfont\itfam=\twelveit
  \textfont\bffam=\fourteenbf \scriptfont\bffam=\twelvebf
        \scriptscriptfont\bffam=\tenbf
  \setbox\strutbox\hbox{\vrule height12.2\p@ depth5\p@ width\z@}%
  \setbox\strutbox@\hbox{\lower.72\normallineskiplimit\vbox{%
        \kern-\normallineskiplimit\copy\strutbox}}%
 \setbox\z@\vbox{\hbox{$($}\kern\z@}\bigsize@=1.7\ht\z@
 \fi
 \normalbaselines\rm\ex@.2326ex\jot4.3\ex@\the\fourteenpoint@}

\def\chapheadfonts{\fourteenpoint\bf}

\font@\seventeenrm=cmr10 scaled\magstep3
\font@\seventeenit=cmti10 scaled\magstep3
\font@\seventeensl=cmsl10 scaled\magstep3
\font@\seventeensmc=cmcsc10 scaled\magstep3
\font@\seventeentt=cmtt10 scaled\magstep3
\font@\seventeenbf=cmbx10 scaled\magstep3
\font@\seventeeni=cmmi10 scaled\magstep3
\font@\seventeensy=cmsy10 scaled\magstep3
\font@\seventeenex=cmex10 scaled\magstep3
\font@\seventeenmsa=msam10 scaled\magstep3
\font@\seventeeneufm=eufm10 scaled\magstep3
\font@\seventeenmsb=msbm10 scaled\magstep3
\newtoks\seventeenpoint@
\def\seventeenpoint{\normalbaselineskip18\p@
 \abovedisplayskip21.6\p@ plus5.2\p@ minus15.4\p@
 \belowdisplayskip\abovedisplayskip
 \abovedisplayshortskip\z@ plus5.2\p@
 \belowdisplayshortskip12.1\p@ plus5.2\p@ minus7\p@
 \textonlyfont@\rm\seventeenrm \textonlyfont@\it\seventeenit
 \textonlyfont@\sl\seventeensl \textonlyfont@\bf\seventeenbf
 \textonlyfont@\smc\seventeensmc \textonlyfont@\tt\seventeentt
 \textonlyfont@\bsmc\seventeenbsmc
 \ifsyntax@ \def\big##1{{\hbox{$\left##1\right.$}}}%
  \let\Big\big \let\bigg\big \let\Bigg\big
 \else
  \textfont\z@=\seventeenrm  \scriptfont\z@=\fourteenrm  \scriptscriptfont\z@=\twelverm
  \textfont\@ne=\seventeeni  \scriptfont\@ne=\fourteeni  \scriptscriptfont\@ne=\twelvei
  \textfont\tw@=\seventeensy \scriptfont\tw@=\fourteensy \scriptscriptfont\tw@=\twelvesy
  \textfont\thr@@=\seventeenex \scriptfont\thr@@=\fourteenex
        \scriptscriptfont\thr@@=\fourteenex
  \textfont\itfam=\seventeenit \scriptfont\itfam=\fourteenit
        \scriptscriptfont\itfam=\fourteenit
  \textfont\bffam=\seventeenbf \scriptfont\bffam=\fourteenbf
        \scriptscriptfont\bffam=\twelvebf
  \setbox\strutbox\hbox{\vrule height14.6\p@ depth6\p@ width\z@}%
  \setbox\strutbox@\hbox{\lower.86\normallineskiplimit\vbox{%
        \kern-\normallineskiplimit\copy\strutbox}}%
 \setbox\z@\vbox{\hbox{$($}\kern\z@}\bigsize@=2\ht\z@
 \fi
 \normalbaselines\rm\ex@.2326ex\jot5.2\ex@\the\seventeenpoint@}

\font@\rrrrrm=cmr10 scaled\magstep4
\font@\bigtitlefont=cmbx10 scaled\magstep4

\parindent1pc
\normallineskiplimit\p@
\newdimen\indenti \indenti=2pc
\def\pageheight#1{\vsize#1}
\def\pagewidth#1{\hsize#1%
   \captionwidth@\hsize \advance\captionwidth@-2\indenti}
\pagewidth{30pc} \pageheight{47pc}
\def\topmatter{%
 \ifx\undefined\msafam
 \else\font@\eightmsa=msam8 \font@\sixmsa=msam6
   \ifsyntax@\else \addto\tenpoint{\textfont\msafam=\tenmsa
              \scriptfont\msafam=\sevenmsa \scriptscriptfont\msafam=\fivemsa}%
     \addto\eightpoint{\textfont\msafam=\eightmsa \scriptfont\msafam=\sixmsa
              \scriptscriptfont\msafam=\fivemsa}%
   \fi
 \fi
 \ifx\undefined\msbfam
 \else\font@\eightmsb=msbm8 \font@\sixmsb=msbm6
   \ifsyntax@\else \addto\tenpoint{\textfont\msbfam=\tenmsb
         \scriptfont\msbfam=\sevenmsb \scriptscriptfont\msbfam=\fivemsb}%
     \addto\eightpoint{\textfont\msbfam=\eightmsb \scriptfont\msbfam=\sixmsb
         \scriptscriptfont\msbfam=\fivemsb}%
   \fi
 \fi
 \ifx\undefined\eufmfam
 \else \font@\eighteufm=eufm8 \font@\sixeufm=eufm6
   \ifsyntax@\else \addto\tenpoint{\textfont\eufmfam=\teneufm
       \scriptfont\eufmfam=\seveneufm \scriptscriptfont\eufmfam=\fiveeufm}%
     \addto\eightpoint{\textfont\eufmfam=\eighteufm
       \scriptfont\eufmfam=\sixeufm \scriptscriptfont\eufmfam=\fiveeufm}%
   \fi
 \fi
 \ifx\undefined\eufbfam
 \else \font@\eighteufb=eufb8 \font@\sixeufb=eufb6
   \ifsyntax@\else \addto\tenpoint{\textfont\eufbfam=\teneufb
      \scriptfont\eufbfam=\seveneufb \scriptscriptfont\eufbfam=\fiveeufb}%
    \addto\eightpoint{\textfont\eufbfam=\eighteufb
      \scriptfont\eufbfam=\sixeufb \scriptscriptfont\eufbfam=\fiveeufb}%
   \fi
 \fi
 \ifx\undefined\eusmfam
 \else \font@\eighteusm=eusm8 \font@\sixeusm=eusm6
   \ifsyntax@\else \addto\tenpoint{\textfont\eusmfam=\teneusm
       \scriptfont\eusmfam=\seveneusm \scriptscriptfont\eusmfam=\fiveeusm}%
     \addto\eightpoint{\textfont\eusmfam=\eighteusm
       \scriptfont\eusmfam=\sixeusm \scriptscriptfont\eusmfam=\fiveeusm}%
   \fi
 \fi
 \ifx\undefined\eusbfam
 \else \font@\eighteusb=eusb8 \font@\sixeusb=eusb6
   \ifsyntax@\else \addto\tenpoint{\textfont\eusbfam=\teneusb
       \scriptfont\eusbfam=\seveneusb \scriptscriptfont\eusbfam=\fiveeusb}%
     \addto\eightpoint{\textfont\eusbfam=\eighteusb
       \scriptfont\eusbfam=\sixeusb \scriptscriptfont\eusbfam=\fiveeusb}%
   \fi
 \fi
 \ifx\undefined\eurmfam
 \else \font@\eighteurm=eurm8 \font@\sixeurm=eurm6
   \ifsyntax@\else \addto\tenpoint{\textfont\eurmfam=\teneurm
       \scriptfont\eurmfam=\seveneurm \scriptscriptfont\eurmfam=\fiveeurm}%
     \addto\eightpoint{\textfont\eurmfam=\eighteurm
       \scriptfont\eurmfam=\sixeurm \scriptscriptfont\eurmfam=\fiveeurm}%
   \fi
 \fi
 \ifx\undefined\eurbfam
 \else \font@\eighteurb=eurb8 \font@\sixeurb=eurb6
   \ifsyntax@\else \addto\tenpoint{\textfont\eurbfam=\teneurb
       \scriptfont\eurbfam=\seveneurb \scriptscriptfont\eurbfam=\fiveeurb}%
    \addto\eightpoint{\textfont\eurbfam=\eighteurb
       \scriptfont\eurbfam=\sixeurb \scriptscriptfont\eurbfam=\fiveeurb}%
   \fi
 \fi
 \ifx\undefined\cmmibfam
 \else \font@\eightcmmib=cmmib8 \font@\sixcmmib=cmmib6
   \ifsyntax@\else \addto\tenpoint{\textfont\cmmibfam=\tencmmib
       \scriptfont\cmmibfam=\sevencmmib \scriptscriptfont\cmmibfam=\fivecmmib}%
    \addto\eightpoint{\textfont\cmmibfam=\eightcmmib
       \scriptfont\cmmibfam=\sixcmmib \scriptscriptfont\cmmibfam=\fivecmmib}%
   \fi
 \fi
 \ifx\undefined\cmbsyfam
 \else \font@\eightcmbsy=cmbsy8 \font@\sixcmbsy=cmbsy6
   \ifsyntax@\else \addto\tenpoint{\textfont\cmbsyfam=\tencmbsy
      \scriptfont\cmbsyfam=\sevencmbsy \scriptscriptfont\cmbsyfam=\fivecmbsy}%
    \addto\eightpoint{\textfont\cmbsyfam=\eightcmbsy
      \scriptfont\cmbsyfam=\sixcmbsy \scriptscriptfont\cmbsyfam=\fivecmbsy}%
   \fi
 \fi
 \let\topmatter\relax}
\def\chapterno@{\uppercase\expandafter{\romannumeral\chaptercount@}}
\newcount\chaptercount@
\def\chapter{\nofrills@{\afterassignment\chapterno@
                        CHAPTER \global\chaptercount@=}\chapter@
 \DNii@##1{\leavevmode\hskip-\leftskip
   \rlap{\vbox to\z@{\vss\centerline{\eightpoint
   \chapter@##1\unskip}\baselineskip2pc\null}}\hskip\leftskip
   \nofrills@false}%
 \FN@\next@}
\newbox\titlebox@

\def\title{\nofrills@{\relax}\title@%
 \DNii@##1\endtitle{\global\setbox\titlebox@\vtop{\tenpoint\bf
 \raggedcenter@\ignorespaces
 \baselineskip1.3\baselineskip\title@{##1}\endgraf}%
 \ifmonograph@ \edef\next{\the\leftheadtoks}\ifx\next\empty
    \leftheadtext{##1}\fi
 \fi
 \edef\next{\the\rightheadtoks}\ifx\next\empty \rightheadtext{##1}\fi
 }\FN@\next@}
\newbox\authorbox@
\def\author#1\endauthor{\global\setbox\authorbox@
 \vbox{\tenpoint\smc\raggedcenter@\ignorespaces
 #1\endgraf}\relaxnext@ \edef\next{\the\leftheadtoks}%
 \ifx\next\empty\leftheadtext{#1}\fi}
\newbox\affilbox@
\def\affil#1\endaffil{\global\setbox\affilbox@
 \vbox{\tenpoint\raggedcenter@\ignorespaces#1\endgraf}}
\newcount\addresscount@
\addresscount@\z@
\def\address#1\endaddress{\global\advance\addresscount@\@ne
  \expandafter\gdef\csname address\number\addresscount@\endcsname
  {\vskip12\p@ minus6\p@\noindent\eightpoint\smc\ignorespaces#1\par}}
\def\email{\nofrills@{\eightpoint{\it E-mail\/}:\enspace}\email@
  \DNii@##1\endemail{%
  \expandafter\gdef\csname email\number\addresscount@\endcsname
  {\def\usualspace{{\it\enspace}}\smallskip\noindent\eightpoint\email@
  \ignorespaces##1\par}}%
 \FN@\next@}
\def\thedate@{}
\def\date#1\enddate{\gdef\thedate@{\tenpoint\ignorespaces#1\unskip}}
\def\thethanks@{}
\def\thanks#1\endthanks{\gdef\thethanks@{\eightpoint\ignorespaces#1.\unskip}}
\def\thekeywords@{}
\def\keywords{\nofrills@{{\it Key words and phrases.\enspace}}\keywords@
 \DNii@##1\endkeywords{\def\thekeywords@{\def\usualspace{{\it\enspace}}%
 \eightpoint\keywords@\ignorespaces##1\unskip.}}%
 \FN@\next@}
\def\thesubjclass@{}
\def\subjclass{\nofrills@{{\rm2010 {\it Mathematics Subject
   Classification\/}.\enspace}}\subjclass@
 \DNii@##1\endsubjclass{\def\thesubjclass@{\def\usualspace
  {{\rm\enspace}}\eightpoint\subjclass@\ignorespaces##1\unskip.}}%
 \FN@\next@}
\newbox\abstractbox@
\def\abstract{\nofrills@{{\smc Abstract.\enspace}}\abstract@
 \DNii@{\setbox\abstractbox@\vbox\bgroup\noindent$$\vbox\bgroup
  \def\envir@{abstract}\advance\hsize-2\indenti
  \usualspace@{{\enspace}}\eightpoint \noindent\abstract@\ignorespaces}%
 \FN@\next@}
\def\endabstract{\par\unskip\egroup$$\egroup}
\def\widestnumber#1#2{\begingroup\let\head\null\let\subhead\empty
   \let\subsubhead\subhead
   \ifx#1\head\global\setbox\tocheadbox@\hbox{#2.\enspace}%
   \else\ifx#1\subhead\global\setbox\tocsubheadbox@\hbox{#2.\enspace}%
   \else\ifx#1\key\bgroup\let\endrefitem@\egroup
        \key#2\endrefitem@\global\refindentwd\wd\keybox@
   \else\ifx#1\no\bgroup\let\endrefitem@\egroup
        \no#2\endrefitem@\global\refindentwd\wd\nobox@
   \else\ifx#1\page\global\setbox\pagesbox@\hbox{\quad\bf#2}%
   \else\ifx#1\item\setboxz@h{#2}\global\rosteritemwd\wdz@
        \global\advance\rosteritemwd by.5\parindent
   \else\message{\string\widestnumber is not defined for this option
   (\string#1)}%
\fi\fi\fi\fi\fi\fi\endgroup}
\newif\ifmonograph@
\def\Monograph{\monograph@true \let\headmark\rightheadtext
  \let\varindent@\indent \def\headfont@{\bf}\def\proclaimheadfont@{\smc}%
  \def\demofont@{\smc}}
\let\varindent@\indent

\newbox\tocheadbox@    \newbox\tocsubheadbox@
\newbox\tocbox@
\def\toc{\toc@{Contents}}
\def\newtocdefs{%
   \def \title##1\endtitle
       {\penaltyandskip@\z@\smallskipamount
        \hangindent\wd\tocheadbox@\noindent{\bf##1}}%
   \def \chapter##1{%
        Chapter \uppercase\expandafter{\romannumeral##1.\unskip}\enspace}%
   \def \specialhead##1\endspecialhead
       {\par\hangindent\wd\tocheadbox@ \noindent##1\par}%
   \def \head##1 ##2\endhead
       {\par\hangindent\wd\tocheadbox@ \noindent
        \if\notempty{##1}\hbox to\wd\tocheadbox@{\hfil##1\enspace}\fi
        ##2\par}%
   \def \subhead##1 ##2\endsubhead
       {\par\vskip-\parskip {\normalbaselines
        \advance\leftskip\wd\tocheadbox@
        \hangindent\wd\tocsubheadbox@ \noindent
        \if\notempty{##1}\hbox to\wd\tocsubheadbox@{##1\unskip\hfil}\fi
         ##2\par}}%
   \def \subsubhead##1 ##2\endsubsubhead
       {\par\vskip-\parskip {\normalbaselines
        \advance\leftskip\wd\tocheadbox@
        \hangindent\wd\tocsubheadbox@ \noindent
        \if\notempty{##1}\hbox to\wd\tocsubheadbox@{##1\unskip\hfil}\fi
        ##2\par}}}
\def\toc@#1{\relaxnext@
   \def\page##1%
       {\unskip\penalty0\null\hfil
        \rlap{\hbox to\wd\pagesbox@{\quad\hfil##1}}\hfilneg\penalty\@M}%
 \DN@{\ifx\next\nofrills\DN@\nofrills{\nextii@}%
      \else\DN@{\nextii@{{#1}}}\fi
      \next@}%
 \DNii@##1{%
\ifmonograph@\bgroup\else\setbox\tocbox@\vbox\bgroup
   \centerline{\headfont@\ignorespaces##1\unskip}\nobreak
   \vskip\belowheadskip \fi
   \setbox\tocheadbox@\hbox{0.\enspace}%
   \setbox\tocsubheadbox@\hbox{0.0.\enspace}%
   \leftskip\indenti \rightskip\leftskip
   \setbox\pagesbox@\hbox{\bf\quad000}\advance\rightskip\wd\pagesbox@
   \newtocdefs
 }%
 \FN@\next@}
\def\endtoc{\par\egroup}
\let\pretitle\relax
\let\preauthor\relax
\let\preaffil\relax
\let\predate\relax
\let\preabstract\relax
\let\prepaper\relax
\def\dedicatory #1\enddedicatory{\def\preabstract{{\medskip
  \eightpoint\it \raggedcenter@#1\endgraf}}}
\def\thetranslator@{}
\def\translator#1\endtranslator{\def\thetranslator@{\nobreak\medskip
 \line{\eightpoint\hfil Translated by \uppercase{#1}\qquad\qquad}\nobreak}}
\outer\def\endtopmatter{\runaway@{abstract}%
 \edef\next{\the\leftheadtoks}\ifx\next\empty
  \expandafter\leftheadtext\expandafter{\the\rightheadtoks}\fi
 \ifmonograph@\else
   \ifx\thesubjclass@\empty\else \makefootnote@{}{\thesubjclass@}\fi
   \ifx\thekeywords@\empty\else \makefootnote@{}{\thekeywords@}\fi
   \ifx\thethanks@\empty\else \makefootnote@{}{\thethanks@}\fi
 \fi
  \pretitle
  \ifmonograph@ \topskip7pc \else \topskip4pc \fi
  \box\titlebox@
  \topskip10pt
  \preauthor
  \ifvoid\authorbox@\else \vskip2.5pc plus1pc \unvbox\authorbox@\fi
  \preaffil
  \ifvoid\affilbox@\else \vskip1pc plus.5pc \unvbox\affilbox@\fi
  \predate
  \ifx\thedate@\empty\else \vskip1pc plus.5pc \line{\hfil\thedate@\hfil}\fi
  \preabstract
  \ifvoid\abstractbox@\else \vskip1.5pc plus.5pc \unvbox\abstractbox@ \fi
  \ifvoid\tocbox@\else\vskip1.5pc plus.5pc \unvbox\tocbox@\fi
  \prepaper
  \vskip2pc plus1pc
}
\def\document{\let\fontlist@\relax\let\alloclist@\relax
  \tenpoint}

\newskip\aboveheadskip       \aboveheadskip1.8\bigskipamount
\newdimen\belowheadskip      \belowheadskip1.8\medskipamount

\def\headfont@{\smc}
\def\penaltyandskip@#1#2{\relax\ifdim\lastskip<#2\relax\removelastskip
      \ifnum#1=\z@\else\penalty@#1\relax\fi\vskip#2%
  \else\ifnum#1=\z@\else\penalty@#1\relax\fi\fi}
\def\nobreak{\penalty\@M
  \ifvmode\def\penalty@{\let\penalty@\penalty\count@@@}%
  \everypar{\let\penalty@\penalty\everypar{}}\fi}
\let\penalty@\penalty
\def\heading#1\endheading{\head#1\endhead}

\def\specialheadfont@{\bf}
\outer\def\specialhead{\par\penaltyandskip@{-200}\aboveheadskip
  \begingroup\interlinepenalty\@M\rightskip\z@ plus\hsize \let\\\linebreak
  \specialheadfont@\noindent\ignorespaces}
\def\endspecialhead{\par\endgroup\nobreak\vskip\belowheadskip}
\let\headmark\eat@
\newskip\subheadskip       \subheadskip\medskipamount
\def\subheadfont@{\bf}
\outer\def\subhead{\nofrills@{.\enspace}\subhead@
 \DNii@##1\endsubhead{\par\penaltyandskip@{-100}\subheadskip
  \varindent@{\usualspace@{{\subheadfont@\enspace}}%
 \subheadfont@\ignorespaces##1\unskip\subhead@}\ignorespaces}%
 \FN@\next@}
\outer\def\subsubhead{\nofrills@{.\enspace}\subsubhead@
 \DNii@##1\endsubsubhead{\par\penaltyandskip@{-50}\medskipamount
      {\usualspace@{{\it\enspace}}%
  \it\ignorespaces##1\unskip\subsubhead@}\ignorespaces}%
 \FN@\next@}
\def\proclaimheadfont@{\bf}
\outer\def\proclaim{\runaway@{proclaim}\def\envir@{proclaim}%
  \nofrills@{.\enspace}\proclaim@
 \DNii@##1{\penaltyandskip@{-100}\medskipamount\varindent@
   \usualspace@{{\proclaimheadfont@\enspace}}\proclaimheadfont@
   \ignorespaces##1\unskip\proclaim@
  \sl\ignorespaces}%
 \FN@\next@}
\outer\def\endproclaim{\let\envir@\relax\par\rm
  \penaltyandskip@{55}\medskipamount}
\def\demoheadfont@{\it}
\def\demo{\runaway@{proclaim}\nofrills@{.\enspace}\demo@
     \DNii@##1{\par\penaltyandskip@\z@\medskipamount
  {\usualspace@{{\demoheadfont@\enspace}}%
  \varindent@\demoheadfont@\ignorespaces##1\unskip\demo@}\rm
  \ignorespaces}\FN@\next@}
\def\enddemo{\par\medskip}
\def\qed{\ifhmode\unskip\nobreak\fi\quad\ifmmode\square\else$\m@th\square$\fi}
\let\remark\demo
\let\endremark\enddemo
\def\definition{\runaway@{proclaim}%
  \nofrills@{.\demoheadfont@\enspace}\definition@
        \DNii@##1{\penaltyandskip@{-100}\medskipamount
        {\usualspace@{{\demoheadfont@\enspace}}%
        \varindent@\demoheadfont@\ignorespaces##1\unskip\definition@}%
        \rm \ignorespaces}\FN@\next@}


\newdimen\rosteritemwd
\newcount\rostercount@
\newif\iffirstitem@
\let\plainitem@\item
\newtoks\everypartoks@
\def\par@{\everypartoks@\expandafter{\the\everypar}\everypar{}}
\def\roster{\edef\leftskip@{\leftskip\the\leftskip}%
 \relaxnext@
 \rostercount@\z@  
 \def\item{\FN@\rosteritem@}%
 \DN@{\ifx\next\runinitem\let\next@\nextii@\else
  \let\next@\nextiii@\fi\next@}%
 \DNii@\runinitem  
  {\unskip  
   \DN@{\ifx\next[\let\next@\nextii@\else
    \ifx\next"\let\next@\nextiii@\else\let\next@\nextiv@\fi\fi\next@}%
   \DNii@[####1]{\rostercount@####1\relax
    \enspace{\rm(\number\rostercount@)}~\ignorespaces}%
   \def\nextiii@"####1"{\enspace{\rm####1}~\ignorespaces}%
   \def\nextiv@{\enspace{\rm(1)}\rostercount@\@ne~}%
   \par@\firstitem@false  
   \FN@\next@}%
 \def\nextiii@{\par\par@  
  \penalty\@m\smallskip\vskip-\parskip
  \firstitem@true}%
 \FN@\next@}
\def\rosteritem@{\iffirstitem@\firstitem@false\else\par\vskip-\parskip\fi
 \leftskip3\parindent\noindent  
 \DNii@[##1]{\rostercount@##1\relax
  \llap{\hbox to2.5\parindent{\hss\rm(\number\rostercount@)}%
   \hskip.5\parindent}\ignorespaces}%
 \def\nextiii@"##1"{%
  \llap{\hbox to2.5\parindent{\hss\rm##1}\hskip.5\parindent}\ignorespaces}%
 \def\nextiv@{\advance\rostercount@\@ne
  \llap{\hbox to2.5\parindent{\hss\rm(\number\rostercount@)}%
   \hskip.5\parindent}}%
 \ifx\next[\let\next@\nextii@\else\ifx\next"\let\next@\nextiii@\else
  \let\next@\nextiv@\fi\fi\next@}

\newif\ifnextRunin@
\def\endroster{\relaxnext@
 \par\leftskip@  
 \penalty-50 \vskip-\parskip\smallskip  
 \DN@{\ifx\next\Runinitem\let\next@\relax
  \else\nextRunin@false\let\item\plainitem@  
   \ifx\next\par 
    \DN@\par{\everypar\expandafter{\the\everypartoks@}}%
   \else  
    \DN@{\noindent\everypar\expandafter{\the\everypartoks@}}%
  \fi\fi\next@}%
 \FN@\next@}
\newcount\rosterhangafter@
\def\Runinitem#1\roster\runinitem{\relaxnext@
 \rostercount@\z@ 
 \def\item{\FN@\rosteritem@}%
 \def\runinitem@{#1}%
 \DN@{\ifx\next[\let\next\nextii@\else\ifx\next"\let\next\nextiii@
  \else\let\next\nextiv@\fi\fi\next}%
 \DNii@[##1]{\rostercount@##1\relax
  \def\item@{{\rm(\number\rostercount@)}}\nextv@}%
 \def\nextiii@"##1"{\def\item@{{\rm##1}}\nextv@}%
 \def\nextiv@{\advance\rostercount@\@ne
  \def\item@{{\rm(\number\rostercount@)}}\nextv@}%
 \def\nextv@{\setbox\z@\vbox  
  {\ifnextRunin@\noindent\fi  
  \runinitem@\unskip\enspace\item@~\par  
  \global\rosterhangafter@\prevgraf}%
  \firstitem@false  
  \ifnextRunin@\else\par\fi  
  \hangafter\rosterhangafter@\hangindent3\parindent
  \ifnextRunin@\noindent\fi  
  \runinitem@\unskip\enspace 
  \item@~\ifnextRunin@\else\par@\fi  
  \nextRunin@true\ignorespaces}%
 \FN@\next@}
\def\footmarkform@#1{$\m@th^{#1}$}
\let\thefootnotemark\footmarkform@
\def\makefootnote@#1#2{\insert\footins
 {\interlinepenalty\interfootnotelinepenalty
 \eightpoint\splittopskip\ht\strutbox\splitmaxdepth\dp\strutbox
 \floatingpenalty\@MM\leftskip\z@\rightskip\z@\spaceskip\z@\xspaceskip\z@
 \leavevmode{#1}\footstrut\ignorespaces#2\unskip\lower\dp\strutbox
 \vbox to\dp\strutbox{}}}
\newcount\footmarkcount@
\footmarkcount@\z@
\def\footnotemark{\let\@sf\empty\relaxnext@
 \ifhmode\edef\@sf{\spacefactor\the\spacefactor}\/\fi
 \DN@{\ifx[\next\let\next@\nextii@\else
  \ifx"\next\let\next@\nextiii@\else
  \let\next@\nextiv@\fi\fi\next@}%
 \DNii@[##1]{\footmarkform@{##1}\@sf}%
 \def\nextiii@"##1"{{##1}\@sf}%
 \def\nextiv@{\iffirstchoice@\global\advance\footmarkcount@\@ne\fi
  \footmarkform@{\number\footmarkcount@}\@sf}%
 \FN@\next@}
\def\footnotetext{\relaxnext@
 \DN@{\ifx[\next\let\next@\nextii@\else
  \ifx"\next\let\next@\nextiii@\else
  \let\next@\nextiv@\fi\fi\next@}%
 \DNii@[##1]##2{\makefootnote@{\footmarkform@{##1}}{##2}}%
 \def\nextiii@"##1"##2{\makefootnote@{##1}{##2}}%
 \def\nextiv@##1{\makefootnote@{\footmarkform@{\number\footmarkcount@}}{##1}}%
 \FN@\next@}
\def\footnote{\let\@sf\empty\relaxnext@
 \ifhmode\edef\@sf{\spacefactor\the\spacefactor}\/\fi
 \DN@{\ifx[\next\let\next@\nextii@\else
  \ifx"\next\let\next@\nextiii@\else
  \let\next@\nextiv@\fi\fi\next@}%
 \DNii@[##1]##2{\footnotemark[##1]\footnotetext[##1]{##2}}%
 \def\nextiii@"##1"##2{\footnotemark"##1"\footnotetext"##1"{##2}}%
 \def\nextiv@##1{\footnotemark\footnotetext{##1}}%
 \FN@\next@}
\def\adjustfootnotemark#1{\advance\footmarkcount@#1\relax}
\def\footnoterule{\kern-3\p@
  \hrule width 5pc\kern 2.6\p@} 
\def\captionfont@{\smc}
\def\topcaption#1#2\endcaption{%
  {\dimen@\hsize \advance\dimen@-\captionwidth@
   \rm\raggedcenter@ \advance\leftskip.5\dimen@ \rightskip\leftskip
  {\captionfont@#1}%
  \if\notempty{#2}.\enspace\ignorespaces#2\fi
  \endgraf}\nobreak\bigskip}
\def\botcaption#1#2\endcaption{%
  \nobreak\bigskip
  \setboxz@h{\captionfont@#1\if\notempty{#2}.\enspace\rm#2\fi}%
  {\dimen@\hsize \advance\dimen@-\captionwidth@
   \leftskip.5\dimen@ \rightskip\leftskip
   \noindent \ifdim\wdz@>\captionwidth@ 
   \else\hfil\fi 
  {\captionfont@#1}\if\notempty{#2}.\enspace\rm#2\fi\endgraf}}
\def\@ins{\par\begingroup\def\vspace##1{\vskip##1\relax}%
  \def\captionwidth##1{\captionwidth@##1\relax}%
  \setbox\z@\vbox\bgroup} 
\def\block{\RIfMIfI@\nondmatherr@\block\fi
       \else\ifvmode\vskip\abovedisplayskip\noindent\fi
        $$\def\endblock{\par\egroup$$}\fi
  \vbox\bgroup\advance\hsize-2\indenti\noindent}
\def\endblock{\par\egroup}
\def\cite#1{{\rm[{\citefont@\m@th#1}]}}
\def\citefont@{\rm}
\def\refsfont@{\eightpoint}
\outer\def\Refs{\runaway@{proclaim}%
 \relaxnext@ \DN@{\ifx\next\nofrills\DN@\nofrills{\nextii@}\else
  \DN@{\nextii@{References}}\fi\next@}%
 \DNii@##1{\penaltyandskip@{-200}\aboveheadskip
  \line{\hfil\headfont@\ignorespaces##1\unskip\hfil}\nobreak
  \vskip\belowheadskip
  \begingroup\refsfont@\sfcode`.=\@m}%
 \FN@\next@}
\def\endRefs{\par\endgroup}
\newbox\nobox@            \newbox\keybox@           \newbox\bybox@
\newbox\paperbox@         \newbox\paperinfobox@     \newbox\jourbox@
\newbox\volbox@           \newbox\issuebox@         \newbox\yrbox@
\newbox\pagesbox@         \newbox\bookbox@          \newbox\bookinfobox@
\newbox\publbox@          \newbox\publaddrbox@      \newbox\finalinfobox@
\newbox\edsbox@           \newbox\langbox@
\newif\iffirstref@        \newif\iflastref@
\newif\ifprevjour@        \newif\ifbook@            \newif\ifprevinbook@
\newif\ifquotes@          \newif\ifbookquotes@      \newif\ifpaperquotes@
\newdimen\bysamerulewd@
\setboxz@h{\refsfont@\kern3em}
\bysamerulewd@\wdz@
\newdimen\refindentwd
\setboxz@h{\refsfont@ 00. }
\refindentwd\wdz@
\outer\def\ref{\begingroup \noindent\hangindent\refindentwd
 \firstref@true \def\nofrills{\def\refkern@{\kern3sp}}%
 \ref@}
\def\ref@{\book@false \bgroup\let\endrefitem@\egroup \ignorespaces}
\def\moreref{\endrefitem@\endref@\firstref@false\ref@}%
\def\transl{\endrefitem@\endref@\firstref@false
  \book@false
  \prepunct@
  \setboxz@h\bgroup \aftergroup\unhbox\aftergroup\z@
    \def\endrefitem@{\unskip\refkern@\egroup}\ignorespaces}%
\def\emptyifempty@{\dimen@\wd\currbox@
  \advance\dimen@-\wd\z@ \advance\dimen@-.1\p@
  \ifdim\dimen@<\z@ \setbox\currbox@\copy\voidb@x \fi}
\let\refkern@\relax
\def\endrefitem@{\unskip\refkern@\egroup
  \setboxz@h{\refkern@}\emptyifempty@}\ignorespaces
\def\refdef@#1#2#3{\edef\next@{\noexpand\endrefitem@
  \let\noexpand\currbox@\csname\expandafter\eat@\string#1box@\endcsname
    \noexpand\setbox\noexpand\currbox@\hbox\bgroup}%
  \toks@\expandafter{\next@}%
  \if\notempty{#2#3}\toks@\expandafter{\the\toks@
  \def\endrefitem@{\unskip#3\refkern@\egroup
  \setboxz@h{#2#3\refkern@}\emptyifempty@}#2}\fi
  \toks@\expandafter{\the\toks@\ignorespaces}%
  \edef#1{\the\toks@}}
\refdef@\no{}{. }
\refdef@\key{[\m@th}{] }
\refdef@\by{}{}
\def\bysame{\by\hbox to\bysamerulewd@{\hrulefill}\thinspace
   \kern0sp}
\def\manyby{\message{\string\manyby is no longer necessary; \string\by
  can be used instead, starting with version 2.0 of \styname.STY}\by}
\refdef@\paper{\ifpaperquotes@``\fi\it}{}
\refdef@\paperinfo{}{}
\def\jour{\endrefitem@\let\currbox@\jourbox@
  \setbox\currbox@\hbox\bgroup
  \def\endrefitem@{\unskip\refkern@\egroup
    \setboxz@h{\refkern@}\emptyifempty@
    \ifvoid\jourbox@\else\prevjour@true\fi}%
\ignorespaces}
\refdef@\vol{\ifbook@\else\bf\fi}{}
\refdef@\issue{no. }{}
\refdef@\yr{}{}
\refdef@\pages{}{}
\def\page{\endrefitem@\def\pp@{\def\pp@{pp.~}p.~}\let\currbox@\pagesbox@
  \setbox\currbox@\hbox\bgroup\ignorespaces}
\def\pp@{pp.~}
\def\book{\endrefitem@ \let\currbox@\bookbox@
 \setbox\currbox@\hbox\bgroup\def\endrefitem@{\unskip\refkern@\egroup
  \setboxz@h{\ifbookquotes@``\fi}\emptyifempty@
  \ifvoid\bookbox@\else\book@true\fi}%
  \ifbookquotes@``\fi\it\ignorespaces}
\def\inbook{\endrefitem@
  \let\currbox@\bookbox@\setbox\currbox@\hbox\bgroup
  \def\endrefitem@{\unskip\refkern@\egroup
  \setboxz@h{\ifbookquotes@``\fi}\emptyifempty@
  \ifvoid\bookbox@\else\book@true\previnbook@true\fi}%
  \ifbookquotes@``\fi\ignorespaces}
\refdef@\eds{(}{, eds.)}
\def\ed{\endrefitem@\let\currbox@\edsbox@
 \setbox\currbox@\hbox\bgroup
 \def\endrefitem@{\unskip, ed.)\refkern@\egroup
  \setboxz@h{(, ed.)}\emptyifempty@}(\ignorespaces}
\refdef@\bookinfo{}{}
\refdef@\publ{}{}
\refdef@\publaddr{}{}
\refdef@\finalinfo{}{}
\refdef@\lang{(}{)}

\let\refdef@\relax 
\def\ppunbox@#1{\ifvoid#1\else\prepunct@\unhbox#1\fi}
\def\nocomma@#1{\ifvoid#1\else\changepunct@3\prepunct@\unhbox#1\fi}
\def\changepunct@#1{\ifnum\lastkern<3 \unkern\kern#1sp\fi}
\def\prepunct@{\count@\lastkern\unkern
  \ifnum\lastpenalty=0
    \let\penalty@\relax
  \else
    \edef\penalty@{\penalty\the\lastpenalty\relax}%
  \fi
  \unpenalty
  \let\refspace@\ \ifcase\count@,
\or;\or.\or 
  \or\let\refspace@\relax
  \else,\fi
  \ifquotes@''\quotes@false\fi \penalty@ \refspace@
}
\def\transferpenalty@#1{\dimen@\lastkern\unkern
  \ifnum\lastpenalty=0\unpenalty\let\penalty@\relax
  \else\edef\penalty@{\penalty\the\lastpenalty\relax}\unpenalty\fi
  #1\penalty@\kern\dimen@}
\def\endref{\endrefitem@\lastref@true\endref@
  \par\endgroup \prevjour@false \previnbook@false }
\def\endref@{%
\iffirstref@
  \ifvoid\nobox@\ifvoid\keybox@\indent\fi
  \else\hbox to\refindentwd{\hss\unhbox\nobox@}\fi
  \ifvoid\keybox@
  \else\ifdim\wd\keybox@>\refindentwd
         \box\keybox@
       \else\hbox to\refindentwd{\unhbox\keybox@\hfil}\fi\fi
  \kern4sp\ppunbox@\bybox@
\fi 
  \ifvoid\paperbox@
  \else\prepunct@\unhbox\paperbox@
    \ifpaperquotes@\quotes@true\fi\fi
  \ppunbox@\paperinfobox@
  \ifvoid\jourbox@
    \ifprevjour@ \nocomma@\volbox@
      \nocomma@\issuebox@
      \ifvoid\yrbox@\else\changepunct@3\prepunct@(\unhbox\yrbox@
        \transferpenalty@)\fi
      \ppunbox@\pagesbox@
    \fi 
  \else \prepunct@\unhbox\jourbox@
    \nocomma@\volbox@
    \nocomma@\issuebox@
    \ifvoid\yrbox@\else\changepunct@3\prepunct@(\unhbox\yrbox@
      \transferpenalty@)\fi
    \ppunbox@\pagesbox@
  \fi 
  \ifbook@\prepunct@\unhbox\bookbox@ \ifbookquotes@\quotes@true\fi \fi
  \nocomma@\edsbox@
  \ppunbox@\bookinfobox@
  \ifbook@\ifvoid\volbox@\else\prepunct@ vol.~\unhbox\volbox@
  \fi\fi
  \ppunbox@\publbox@ \ppunbox@\publaddrbox@
  \ifbook@ \ppunbox@\yrbox@
    \ifvoid\pagesbox@
    \else\prepunct@\pp@\unhbox\pagesbox@\fi
  \else
    \ifprevinbook@ \ppunbox@\yrbox@
      \ifvoid\pagesbox@\else\prepunct@\pp@\unhbox\pagesbox@\fi
    \fi \fi
  \ppunbox@\finalinfobox@
  \iflastref@
    \ifvoid\langbox@.\ifquotes@''\fi
    \else\changepunct@2\prepunct@\unhbox\langbox@\fi
  \else
    \ifvoid\langbox@\changepunct@1%
    \else\changepunct@3\prepunct@\unhbox\langbox@
      \changepunct@1\fi
  \fi
}
\outer\def\enddocument{%
 \runaway@{proclaim}%
\ifmonograph@ 
\else
 \nobreak
 \thetranslator@
 \count@\z@ \loop\ifnum\count@<\addresscount@\advance\count@\@ne
 \csname address\number\count@\endcsname
 \csname email\number\count@\endcsname
 \repeat
\fi
 \vfill\supereject\end}

\def\headfont@{\headfonts}
\def\proclaimheadfont@{\bf}
\def\specialheadfont@{\bf}
\def\subheadfont@{\bf}
\def\demoheadfont@{\smc}

\newif\ifThisToToc \ThisToTocfalse
\newif\iftocloaded \tocloadedfalse

\def\C@L{\noexpand\Cal}\def\B@B{\noexpand\Bbb}\def\fR@K{\noexpand\frak}
\def\S@{\noexpand\S}\def\P@P{\noexpand\"}
\def\xpar{\\}

\def\writetoc#1{\iftocloaded\ifThisToToc\begingroup\def\totoc{}
  \def\Cal{\noexpand\C@L}\def\Bbb{\noexpand\B@B}
  \def\frak{\noexpand\fR@K}\def\goth{\frak}\def\S{\noexpand\S@}
  \def\"{\noexpand\P@P}
  \def\xpar{\par\penalty100000 }\def\idx##1{##1}\def\\{\xpar}
  \edef\next@{\write\toc{\noindent#1\leaderfill\noexpand\folio\par}}%
  \next@\endgroup\global\ThisToTocfalse\fi\fi}
\def\leaderfill{\leaders\hbox to 1em{\hss.\hss}\hfill}

\newif\ifindexloaded \indexloadedfalse
\def\idx#1{\ifindexloaded\begingroup\def\ign{}\def\it{}\def\/{}%
 \def\smc{}\def\bf{}\def\tt{}%
 \def\Cal{\noexpand\C@L}\def\Bbb{\noexpand\B@B}%
 \def\frak{\noexpand\fR@K}\def\goth{\frak}\def\S{\noexpand\S@}%
  \def\"{\noexpand\P@P}%
 {\edef\next@{\write\index{#1, \noexpand\folio}}\next@}%
 \endgroup\fi{#1}}
\def\ign#1{}

\def\totoc{\global\ThisToToctrue}

\outer\def\head#1\endhead{\par\penaltyandskip@{-200}\aboveheadskip
 {\headfont@\raggedcenter@\interlinepenalty\@M
 \ignorespaces#1\endgraf}\nobreak
 \vskip\belowheadskip
 \headmark{#1}\writetoc{#1}}

\outer\def\chaphead#1\endchaphead{\par\penaltyandskip@{-200}\aboveheadskip
 {\chapheadfonts\raggedcenter@\interlinepenalty\@M
 \ignorespaces#1\endgraf}\nobreak
 \vskip3\belowheadskip
 \headmark{#1}\writetoc{#1}}

\def\folio{{\foliofont@\ifnum\pageno<\z@ \romannumeral-\pageno
 \else\number\pageno \fi}}
\newtoks\leftheadtoks
\newtoks\rightheadtoks

\def\leftheadtext{\nofrills@{\relax}\lht@
  \DNii@##1{\leftheadtoks\expandafter{\lht@{##1}}%
    \mark{\the\leftheadtoks\noexpand\else\the\rightheadtoks}
    \ifsyntax@\setboxz@h{\def\\{\unskip\space\ignorespaces}%
        \headlinefont@##1}\fi}%
  \FN@\next@}
\def\rightheadtext{\nofrills@{\relax}\rht@
  \DNii@##1{\rightheadtoks\expandafter{\rht@{##1}}%
    \mark{\the\leftheadtoks\noexpand\else\the\rightheadtoks}%
    \ifsyntax@\setboxz@h{\def\\{\unskip\space\ignorespaces}%
        \headlinefont@##1}\fi}%
  \FN@\next@}
\def\NoRunningHeads{\global\runheads@false\global\let\headmark\eat@}

\newif\iffirstpage@     \firstpage@true
\newif\ifrunheads@      \runheads@true

\newdimen\fullhsize \fullhsize=\hsize
\newdimen\fullvsize \fullvsize=\vsize
\def\fullline{\hbox to\fullhsize}

\def\pagenumbers{\gdef\folio{\folio@}}

\let\norunningheads\NoRunningHeads
\def\userunningheads{\global\runheads@true}
\norunningheads

\headline={\def\chapter#1{\chapterno@. }%
  \def\\{\unskip\space\ignorespaces}\ifrunheads@\headlinefont@
    \ifodd\pageno\rightheadline \else\leftheadline\fi
   \else\hfil\fi\ifNoRunHeadline\global\NoRunHeadlinefalse\fi}
\let\folio@\folio
\def\foliofont@{\foliofont}
\def\foliofont{\eightrm}
\def\headlinefont@{\headlinefont}
\def\headlinefont{\eightpoint\smc}
\def\leftheadline{\rlap{\folio}\hfill
   \ifNoRunHeadline\else\iftrue\topmark\fi\fi \hfill}
\def\rightheadline{\hfill\ifNoRunHeadline
   \else \expandafter\fi
  \hfill \llap{\folio}}
\footline={{\eightpoint\bottremark}%
   \ifrunheads@\else\hfil{\let\foliofont\tenrm\folio}\fi\hfil}
\def\bottremark{}
 
\newif\ifNoRunHeadline      
\def\norunninghead{\global\NoRunHeadlinetrue}
\norunninghead

\output={\output@}
%
\newif\ifoffset\offsetfalse
\output={\output@}
\def\output@{%
 \ifoffset 
  \ifodd\count0\advance\hoffset by0.5truecm
   \else\advance\hoffset by-0.5truecm\fi\fi
 \shipout\vbox{%
  \makeheadline \pagebody \makefootline }%
 \advancepageno \ifnum\outputpenalty>-\@MM\else\dosupereject\fi}

\def\indexoutput#1{%
  \ifoffset 
   \ifodd\count0\advance\hoffset by0.5truecm
    \else\advance\hoffset by-0.5truecm\fi\fi
  \shipout\vbox{\makeheadline
  \vbox to\fullvsize{\boxmaxdepth\maxdepth%
  \ifvoid\topins\else\unvbox\topins\fi%
  #1 %
  \ifvoid\footins\else 
    \vskip\skip\footins
    \footnoterule
    \unvbox\footins\fi
  \ifr@ggedbottom \kern-\dimen@ \vfil \fi}%
  \baselineskip2pc
  \makefootline}%
 \global\advance\pageno\@ne
 \ifnum\outputpenalty>-\@MM\else\dosupereject\fi}
 
 \newbox\partialpage \newdimen\halfsize \halfsize=0.5\fullhsize
 \advance\halfsize by-0.5em

 \def\begindoublecolumns{\output={\indexoutput{\unvbox255}}%
   \begingroup \def\line{\fullline}
   \output={\global\setbox\partialpage=\vbox{\unvbox255\bigskip}}\eject
   \output={\doublecolumnout}\hsize=\halfsize \vsize=2\fullvsize}
 \def\enddoublecolumns{\output={\balancecolumns}\eject
  \endgroup \pagegoal=\fullvsize%
  \output={\output@}}
\def\doublecolumnout{\splittopskip=\topskip \splitmaxdepth=\maxdepth
  \dimen@=\fullvsize \advance\dimen@ by-\ht\partialpage
  \setbox0=\vsplit255 to \dimen@ \setbox2=\vsplit255 to \dimen@
  \indexoutput{\pagesofar} \unvbox255 \penalty\outputpenalty}
\def\pagesofar{\unvbox\partialpage
  \wd0=\hsize \wd2=\hsize \hbox to\fullhsize{\box0\hfil\box2}}
\def\balancecolumns{\setbox0=\vbox{\unvbox255} \dimen@=\ht0
  \advance\dimen@ by\topskip \advance\dimen@ by-\baselineskip
  \divide\dimen@ by2 \splittopskip=\topskip
  {\vbadness=10000 \loop \global\setbox3=\copy0
    \global\setbox1=\vsplit3 to\dimen@
    \ifdim\ht3>\dimen@ \global\advance\dimen@ by1pt \repeat}
  \setbox0=\vbox to\dimen@{\unvbox1} \setbox2=\vbox to\dimen@{\unvbox3}
  \pagesofar}

\tenpoint
\catcode`\@=\active

\def\smallheadings{\let\chapheadfonts\tenpoint\let\headfonts\tenpoint}

\tenpoint
\catcode`\@=\active

\def\LL{\leavevmode\setbox0=\hbox{L}\hbox to\wd0{\hss\char'40L}}
\def\al{\alpha}

\def\la{\lambda}


\def\today{\ifcase\month\or
 January\or February\or March\or April\or May\or June\or
 July\or August\or September\or October\or November\or December\fi
 \space\number\day, \number\year}

\def\({\left(}
\def\){\right)}
\def\[{\left[}
\def\]{\right]}

\def\maj{\operatorname{maj}}

\def\3{\ss}
\catcode`\@=11
\def\dddot#1{\vbox{\ialign{##\crcr
      .\hskip-.5pt.\hskip-.5pt.\crcr\noalign{\kern1.5\p@\nointerlineskip}
      $\hfil\displaystyle{#1}\hfil$\crcr}}}

\newif\iftab@\tab@false
\newif\ifvtab@\vtab@false
\def\tab{\bgroup\tab@true\vtab@false\vst@bfalse\Strich@false%
   \def\\{\global\hline@@false%
     \ifhline@\global\hline@false\global\hline@@true\fi\cr}
   \edef\l@{\the\leftskip}\ialign\bgroup\hskip\l@##\hfil&&##\hfil\cr}
\def\endtab{\cr\egroup\egroup}
\def\vtab{\vtop\bgroup\vst@bfalse\vtab@true\tab@true\Strich@false%
   \bgroup\def\\{\cr}\ialign\bgroup&##\hfil\cr}
\def\endvtab{\cr\egroup\egroup\egroup}
\def\stab{\D@cke0.5pt\null 
 \bgroup\tab@true\vtab@false\vst@bfalse\Strich@true\Let@@\vspace@
 \normalbaselines\offinterlineskip
  \openup\spreadmlines@
 \edef\l@{\the\leftskip}\ialign
 \bgroup\hskip\l@##\hfil&&##\hfil\crcr}
\def\endstab{\crcr\egroup
 \egroup}
\newif\ifvst@b\vst@bfalse
\def\vstab{\D@cke0.5pt\null
 \vtop\bgroup\tab@true\vtab@false\vst@btrue\Strich@true\bgroup\Let@@\vspace@
 \normalbaselines\offinterlineskip
  \openup\spreadmlines@\bgroup}
\def\endvstab{\crcr\egroup\egroup
 \egroup\tab@false\Strich@false}

\newdimen\htstrut@
\htstrut@8.5\p@
\newdimen\htStrut@
\htStrut@12\p@
\newdimen\dpstrut@
\dpstrut@3.5\p@
\newdimen\dpStrut@
\dpStrut@3.5\p@
\def\openup{\afterassignment\@penup\dimen@=}
\def\@penup{\advance\lineskip\dimen@
  \advance\baselineskip\dimen@
  \advance\lineskiplimit\dimen@
  \divide\dimen@ by2
  \advance\htstrut@\dimen@
  \advance\htStrut@\dimen@
  \advance\dpstrut@\dimen@
  \advance\dpStrut@\dimen@}
\def\Let@@{\relax%
    \def\\{\global\hline@@false%
     \ifhline@\global\hline@false\global\hline@@true\fi\cr}%
    \iffalse}\fi}
\def\matrix{\null\,\vcenter\bgroup
 \tab@false\vtab@false\vst@bfalse\Strich@false\Let@@\vspace@
 \normalbaselines\openup\spreadmlines@\ialign
 \bgroup\hfil$\m@th##$\hfil&&\quad\hfil$\m@th##$\hfil\crcr
 \Mathstrut@\crcr\noalign{\kern-\baselineskip}}
\def\endmatrix{\crcr\Mathstrut@\crcr\noalign{\kern-\baselineskip}\egroup
 \egroup\,}
\def\smatrix{\D@cke0.5pt\null\,
 \vcenter\bgroup\tab@false\vtab@false\vst@bfalse\Strich@true\Let@@\vspace@
 \normalbaselines\offinterlineskip
  \openup\spreadmlines@\ialign
 \bgroup\hfil$\m@th##$\hfil&&\quad\hfil$\m@th##$\hfil\crcr}
\def\endsmatrix{\crcr\egroup
 \egroup\,\Strich@false}
\newdimen\D@cke
\def\Dicke#1{\global\D@cke#1}
\newtoks\tabs@\tabs@{&}
\newif\ifStrich@\Strich@false
\newif\iff@rst

\def\Stricherr@{\iftab@\ifvtab@\errmessage{\noexpand\s not allowed
     here. Use \noexpand\vstab!}%
  \else\errmessage{\noexpand\s not allowed here. Use \noexpand\stab!}%
  \fi\else\errmessage{\noexpand\s not allowed
     here. Use \noexpand\smatrix!}\fi}
\def\format{\ifvst@b\else\crcr\fi\egroup\iffalse{\fi\ifnum`}=0 \fi\format@}
\def\format@#1\\{\def\preamble@{#1}%
 \def\Str@chfehlt##1{\ifx##1\s\Stricherr@\fi\ifx##1\\\let\Next\relax%
   \else\let\Next\Str@chfehlt\fi\Next}%
 \def\c{\hfil\noexpand\ifhline@@\hbox{\vrule height\htStrut@%
   depth\dpstrut@ width\z@}\noexpand\fi%
   \ifStrich@\hbox{\vrule height\htstrut@ depth\dpstrut@ width\z@}%
   \fi\iftab@\else$\m@th\fi\the\hashtoks@\iftab@\else$\fi\hfil}%
 \def\r{\hfil\noexpand\ifhline@@\hbox{\vrule height\htStrut@%
   depth\dpstrut@ width\z@}\noexpand\fi%
   \ifStrich@\hbox{\vrule height\htstrut@ depth\dpstrut@ width\z@}%
   \fi\iftab@\else$\m@th\fi\the\hashtoks@\iftab@\else$\fi}%
 \def\l{\noexpand\ifhline@@\hbox{\vrule height\htStrut@%
   depth\dpstrut@ width\z@}\noexpand\fi%
   \ifStrich@\hbox{\vrule height\htstrut@ depth\dpstrut@ width\z@}%
   \fi\iftab@\else$\m@th\fi\the\hashtoks@\iftab@\else$\fi\hfil}%
 \def\s{\ifStrich@\ \the\tabs@\vrule width\D@cke\the\hashtoks@%
          \fi\the\tabs@\ }%
 \def\sa{\ifStrich@\vrule width\D@cke\the\hashtoks@%
            \the\tabs@\ %
            \fi}%
 \def\se{\ifStrich@\ \the\tabs@\vrule width\D@cke\the\hashtoks@\fi}%
 \def\cd{\hfil\noexpand\ifhline@@\hbox{\vrule height\htStrut@%
   depth\dpstrut@ width\z@}\noexpand\fi%
   \ifStrich@\hbox{\vrule height\htstrut@ depth\dpstrut@ width\z@}%
   \fi$\dsize\m@th\the\hashtoks@$\hfil}%
 \def\rd{\hfil\noexpand\ifhline@@\hbox{\vrule height\htStrut@%
   depth\dpstrut@ width\z@}\noexpand\fi%
   \ifStrich@\hbox{\vrule height\htstrut@ depth\dpstrut@ width\z@}%
   \fi$\dsize\m@th\the\hashtoks@$}%
 \def\ld{\noexpand\ifhline@@\hbox{\vrule height\htStrut@%
   depth\dpstrut@ width\z@}\noexpand\fi%
   \ifStrich@\hbox{\vrule height\htstrut@ depth\dpstrut@ width\z@}%
   \fi$\dsize\m@th\the\hashtoks@$\hfil}%
 \ifStrich@\else\Str@chfehlt#1\\\fi%
 \setbox\z@\hbox{\xdef\Preamble@{\preamble@}}\ifnum`{=0 \fi\iffalse}\fi
 \ialign\bgroup\span\Preamble@\crcr}
\newif\ifhline@\hline@false
\newif\ifhline@@\hline@@false
\def\hlinefor#1{\multispan@{\strip@#1 }\leaders\hrule height\D@cke\hfill%
    \global\hline@true\ignorespaces}
\def\Item "#1"{\par\noindent\hangindent2\parindent%
  \hangafter1\setbox0\hbox{\rm#1\enspace}\ifdim\wd0>2\parindent%
  \box0\else\hbox to 2\parindent{\rm#1\hfil}\fi\ignorespaces}
\def\ITEM #1"#2"{\par\noindent\hangafter1\hangindent#1%
  \setbox0\hbox{\rm#2\enspace}\ifdim\wd0>#1%
  \box0\else\hbox to 0pt{\rm#2\hss}\hskip#1\fi\ignorespaces}
\def\item"#1"{\par\noindent\hang%
  \setbox0=\hbox{\rm#1\enspace}\ifdim\wd0>\the\parindent%
  \box0\else\hbox to \parindent{\rm#1\hfil}\enspace\fi\ignorespaces}
\let\plainitem@\item
\catcode`\@=13

\magnification1200

\catcode`\@=11
\font\tenln    = line10
\font\tenlnw   = linew10

\newskip\Einheit \Einheit=0.5cm
\newcount\xcoord \newcount\ycoord
\newdimen\xdim \newdimen\ydim \newdimen\PfadD@cke \newdimen\Pfadd@cke

\newcount\@tempcnta
\newcount\@tempcntb

\newdimen\@tempdima
\newdimen\@tempdimb

\newdimen\@wholewidth
\newdimen\@halfwidth

\newcount\@xarg
\newcount\@yarg
\newcount\@yyarg
\newbox\@linechar
\newbox\@tempboxa
\newdimen\@linelen
\newdimen\@clnwd
\newdimen\@clnht

\newif\if@negarg

\def\@whilenoop#1{}
\def\@whiledim#1\do #2{\ifdim #1\relax#2\@iwhiledim{#1\relax#2}\fi}
\def\@iwhiledim#1{\ifdim #1\let\@nextwhile=\@iwhiledim
        \else\let\@nextwhile=\@whilenoop\fi\@nextwhile{#1}}

\def\@whileswnoop#1\fi{}
\def\@whilesw#1\fi#2{#1#2\@iwhilesw{#1#2}\fi\fi}
\def\@iwhilesw#1\fi{#1\let\@nextwhile=\@iwhilesw
         \else\let\@nextwhile=\@whileswnoop\fi\@nextwhile{#1}\fi}

\def\thinlines{\let\@linefnt\tenln \let\@circlefnt\tencirc
  \@wholewidth\fontdimen8\tenln \@halfwidth .5\@wholewidth}
\def\thicklines{\let\@linefnt\tenlnw \let\@circlefnt\tencircw
  \@wholewidth\fontdimen8\tenlnw \@halfwidth .5\@wholewidth}
\thinlines

\PfadD@cke1pt \Pfadd@cke0.5pt
\def\PfadDicke#1{\PfadD@cke#1 \divide\PfadD@cke by2 \Pfadd@cke\PfadD@cke \multiply\PfadD@cke by2}
\long\def\LOOP#1\REPEAT{\def\BODY{#1}\ITERATE}
\def\ITERATE{\BODY \let\next\ITERATE \else\let\next\relax\fi \next}
\let\REPEAT=\fi
\def\Punkt{\hbox{\raise-2pt\hbox to0pt{\hss$\ssize\bullet$\hss}}}
\def\DuennPunkt(#1,#2){\unskip
  \raise#2 \Einheit\hbox to0pt{\hskip#1 \Einheit
          \raise-2.5pt\hbox to0pt{\hss$\bullet$\hss}\hss}}
\def\NormalPunkt(#1,#2){\unskip
  \raise#2 \Einheit\hbox to0pt{\hskip#1 \Einheit
          \raise-3pt\hbox to0pt{\hss\twelvepoint$\bullet$\hss}\hss}}
\def\DickPunkt(#1,#2){\unskip
  \raise#2 \Einheit\hbox to0pt{\hskip#1 \Einheit
          \raise-4pt\hbox to0pt{\hss\fourteenpoint$\bullet$\hss}\hss}}
\def\Kreis(#1,#2){\unskip
  \raise#2 \Einheit\hbox to0pt{\hskip#1 \Einheit
          \raise-4pt\hbox to0pt{\hss\fourteenpoint$\circ$\hss}\hss}}

\def\Line@(#1,#2)#3{\@xarg #1\relax \@yarg #2\relax
\@linelen=#3\Einheit
\ifnum\@xarg =0 \@vline
  \else \ifnum\@yarg =0 \@hline \else \@sline\fi
\fi}

\def\@sline{\ifnum\@xarg< 0 \@negargtrue \@xarg -\@xarg \@yyarg -\@yarg
  \else \@negargfalse \@yyarg \@yarg \fi
\ifnum \@yyarg >0 \@tempcnta\@yyarg \else \@tempcnta -\@yyarg \fi
\ifnum\@tempcnta>6 \@badlinearg\@tempcnta0 \fi
\ifnum\@xarg>6 \@badlinearg\@xarg 1 \fi
\setbox\@linechar\hbox{\@linefnt\@getlinechar(\@xarg,\@yyarg)}%
\ifnum \@yarg >0 \let\@upordown\raise \@clnht\z@
   \else\let\@upordown\lower \@clnht \ht\@linechar\fi
\@clnwd=\wd\@linechar
\if@negarg \hskip -\wd\@linechar \def\@tempa{\hskip -2\wd\@linechar}\else
     \let\@tempa\relax \fi
\@whiledim \@clnwd <\@linelen \do
  {\@upordown\@clnht\copy\@linechar
   \@tempa
   \advance\@clnht \ht\@linechar
   \advance\@clnwd \wd\@linechar}%
\advance\@clnht -\ht\@linechar
\advance\@clnwd -\wd\@linechar
\@tempdima\@linelen\advance\@tempdima -\@clnwd
\@tempdimb\@tempdima\advance\@tempdimb -\wd\@linechar
\if@negarg \hskip -\@tempdimb \else \hskip \@tempdimb \fi
\multiply\@tempdima \@m
\@tempcnta \@tempdima \@tempdima \wd\@linechar \divide\@tempcnta \@tempdima
\@tempdima \ht\@linechar \multiply\@tempdima \@tempcnta
\divide\@tempdima \@m
\advance\@clnht \@tempdima
\ifdim \@linelen <\wd\@linechar
   \hskip \wd\@linechar
  \else\@upordown\@clnht\copy\@linechar\fi}

\def\@hline{\ifnum \@xarg <0 \hskip -\@linelen \fi
\vrule height\Pfadd@cke width \@linelen depth\Pfadd@cke
\ifnum \@xarg <0 \hskip -\@linelen \fi}

\def\@getlinechar(#1,#2){\@tempcnta#1\relax\multiply\@tempcnta 8
\advance\@tempcnta -9 \ifnum #2>0 \advance\@tempcnta #2\relax\else
\advance\@tempcnta -#2\relax\advance\@tempcnta 64 \fi
\char\@tempcnta}

\def\Vektor(#1,#2)#3(#4,#5){\unskip\leavevmode
  \xcoord#4\relax \ycoord#5\relax
      \raise\ycoord \Einheit\hbox to0pt{\hskip\xcoord \Einheit
         \Vector@(#1,#2){#3}\hss}}

\def\Vector@(#1,#2)#3{\@xarg #1\relax \@yarg #2\relax
\@tempcnta \ifnum\@xarg<0 -\@xarg\else\@xarg\fi
\ifnum\@tempcnta<5\relax
\@linelen=#3\Einheit
\ifnum\@xarg =0 \@vvector
  \else \ifnum\@yarg =0 \@hvector \else \@svector\fi
\fi
\else\@badlinearg\fi}

\def\@hvector{\@hline\hbox to 0pt{\@linefnt
\ifnum \@xarg <0 \@getlarrow(1,0)\hss\else
    \hss\@getrarrow(1,0)\fi}}

\def\@vvector{\ifnum \@yarg <0 \@downvector \else \@upvector \fi}

\def\@svector{\@sline
\@tempcnta\@yarg \ifnum\@tempcnta <0 \@tempcnta=-\@tempcnta\fi
\ifnum\@tempcnta <5
  \hskip -\wd\@linechar
  \@upordown\@clnht \hbox{\@linefnt  \if@negarg
  \@getlarrow(\@xarg,\@yyarg) \else \@getrarrow(\@xarg,\@yyarg) \fi}%
\else\@badlinearg\fi}

\def\@upline{\hbox to \z@{\hskip -.5\Pfadd@cke \vrule width \Pfadd@cke
   height \@linelen depth \z@\hss}}

\def\@downline{\hbox to \z@{\hskip -.5\Pfadd@cke \vrule width \Pfadd@cke
   height \z@ depth \@linelen \hss}}

\def\@upvector{\@upline\setbox\@tempboxa\hbox{\@linefnt\char'66}\raise
     \@linelen \hbox to\z@{\lower \ht\@tempboxa\box\@tempboxa\hss}}

\def\@downvector{\@downline\lower \@linelen
      \hbox to \z@{\@linefnt\char'77\hss}}

\def\@getlarrow(#1,#2){\ifnum #2 =\z@ \@tempcnta='33\else
\@tempcnta=#1\relax\multiply\@tempcnta \sixt@@n \advance\@tempcnta
-9 \@tempcntb=#2\relax\multiply\@tempcntb \tw@
\ifnum \@tempcntb >0 \advance\@tempcnta \@tempcntb\relax
\else\advance\@tempcnta -\@tempcntb\advance\@tempcnta 64
\fi\fi\char\@tempcnta}

\def\@getrarrow(#1,#2){\@tempcntb=#2\relax
\ifnum\@tempcntb < 0 \@tempcntb=-\@tempcntb\relax\fi
\ifcase \@tempcntb\relax \@tempcnta='55 \or
\ifnum #1<3 \@tempcnta=#1\relax\multiply\@tempcnta
24 \advance\@tempcnta -6 \else \ifnum #1=3 \@tempcnta=49
\else\@tempcnta=58 \fi\fi\or
\ifnum #1<3 \@tempcnta=#1\relax\multiply\@tempcnta
24 \advance\@tempcnta -3 \else \@tempcnta=51\fi\or
\@tempcnta=#1\relax\multiply\@tempcnta
\sixt@@n \advance\@tempcnta -\tw@ \else
\@tempcnta=#1\relax\multiply\@tempcnta
\sixt@@n \advance\@tempcnta 7 \fi\ifnum #2<0 \advance\@tempcnta 64 \fi
\char\@tempcnta}

\def\Diagonale(#1,#2)#3{\unskip\leavevmode
  \xcoord#1\relax \ycoord#2\relax
      \raise\ycoord \Einheit\hbox to0pt{\hskip\xcoord \Einheit
         \Line@(1,1){#3}\hss}}
\def\AntiDiagonale(#1,#2)#3{\unskip\leavevmode
  \xcoord#1\relax \ycoord#2\relax 
      \raise\ycoord \Einheit\hbox to0pt{\hskip\xcoord \Einheit
         \Line@(1,-1){#3}\hss}}
\def\Pfad(#1,#2),#3\endPfad{\unskip\leavevmode
  \xcoord#1 \ycoord#2 \thicklines\ZeichnePfad#3\endPfad\thinlines}
\def\ZeichnePfad#1{\ifx#1\endPfad\let\next\relax
  \else\let\next\ZeichnePfad
    \ifnum#1=1
      \raise\ycoord \Einheit\hbox to0pt{\hskip\xcoord \Einheit
         \vrule height\Pfadd@cke width1 \Einheit depth\Pfadd@cke\hss}%
      \advance\xcoord by 1
    \else\ifnum#1=2
      \raise\ycoord \Einheit\hbox to0pt{\hskip\xcoord \Einheit
        \hbox{\hskip-\PfadD@cke\vrule height1 \Einheit width\PfadD@cke depth0pt}\hss}%
      \advance\ycoord by 1
    \else\ifnum#1=3
      \raise\ycoord \Einheit\hbox to0pt{\hskip\xcoord \Einheit
         \Line@(1,1){1}\hss}
      \advance\xcoord by 1
      \advance\ycoord by 1
    \else\ifnum#1=4
      \raise\ycoord \Einheit\hbox to0pt{\hskip\xcoord \Einheit
         \Line@(1,-1){1}\hss}
      \advance\xcoord by 1
      \advance\ycoord by -1
    \else\ifnum#1=5
      \advance\xcoord by -1
      \raise\ycoord \Einheit\hbox to0pt{\hskip\xcoord \Einheit
         \vrule height\Pfadd@cke width1 \Einheit depth\Pfadd@cke\hss}%
    \else\ifnum#1=6
      \advance\ycoord by -1
      \raise\ycoord \Einheit\hbox to0pt{\hskip\xcoord \Einheit
        \hbox{\hskip-\PfadD@cke\vrule height1 \Einheit width\PfadD@cke depth0pt}\hss}%
    \else\ifnum#1=7
      \advance\xcoord by -1
      \advance\ycoord by -1
      \raise\ycoord \Einheit\hbox to0pt{\hskip\xcoord \Einheit
         \Line@(1,1){1}\hss}
    \else\ifnum#1=8
      \advance\xcoord by -1
      \advance\ycoord by +1
      \raise\ycoord \Einheit\hbox to0pt{\hskip\xcoord \Einheit
         \Line@(1,-1){1}\hss}
    \fi\fi\fi\fi
    \fi\fi\fi\fi
  \fi\next}
\def\hSSchritt{\leavevmode\raise-.4pt\hbox to0pt{\hss.\hss}\hskip.2\Einheit
  \raise-.4pt\hbox to0pt{\hss.\hss}\hskip.2\Einheit
  \raise-.4pt\hbox to0pt{\hss.\hss}\hskip.2\Einheit
  \raise-.4pt\hbox to0pt{\hss.\hss}\hskip.2\Einheit
  \raise-.4pt\hbox to0pt{\hss.\hss}\hskip.2\Einheit}
\def\vSSchritt{\vbox{\baselineskip.2\Einheit\lineskiplimit0pt
\hbox{.}\hbox{.}\hbox{.}\hbox{.}\hbox{.}}}
\def\DSSchritt{\leavevmode\raise-.4pt\hbox to0pt{%
  \hbox to0pt{\hss.\hss}\hskip.2\Einheit
  \raise.2\Einheit\hbox to0pt{\hss.\hss}\hskip.2\Einheit
  \raise.4\Einheit\hbox to0pt{\hss.\hss}\hskip.2\Einheit
  \raise.6\Einheit\hbox to0pt{\hss.\hss}\hskip.2\Einheit
  \raise.8\Einheit\hbox to0pt{\hss.\hss}\hss}}
\def\dSSchritt{\leavevmode\raise-.4pt\hbox to0pt{%
  \hbox to0pt{\hss.\hss}\hskip.2\Einheit
  \raise-.2\Einheit\hbox to0pt{\hss.\hss}\hskip.2\Einheit
  \raise-.4\Einheit\hbox to0pt{\hss.\hss}\hskip.2\Einheit
  \raise-.6\Einheit\hbox to0pt{\hss.\hss}\hskip.2\Einheit
  \raise-.8\Einheit\hbox to0pt{\hss.\hss}\hss}}
\def\SPfad(#1,#2),#3\endSPfad{\unskip\leavevmode
  \xcoord#1 \ycoord#2 \ZeichneSPfad#3\endSPfad}
\def\ZeichneSPfad#1{\ifx#1\endSPfad\let\next\relax
  \else\let\next\ZeichneSPfad
    \ifnum#1=1
      \raise\ycoord \Einheit\hbox to0pt{\hskip\xcoord \Einheit
         \hSSchritt\hss}%
      \advance\xcoord by 1
    \else\ifnum#1=2
      \raise\ycoord \Einheit\hbox to0pt{\hskip\xcoord \Einheit
        \hbox{\hskip-2pt \vSSchritt}\hss}%
      \advance\ycoord by 1
    \else\ifnum#1=3
      \raise\ycoord \Einheit\hbox to0pt{\hskip\xcoord \Einheit
         \DSSchritt\hss}
      \advance\xcoord by 1
      \advance\ycoord by 1
    \else\ifnum#1=4
      \raise\ycoord \Einheit\hbox to0pt{\hskip\xcoord \Einheit
         \dSSchritt\hss}
      \advance\xcoord by 1
      \advance\ycoord by -1
    \else\ifnum#1=5
      \advance\xcoord by -1
      \raise\ycoord \Einheit\hbox to0pt{\hskip\xcoord \Einheit
         \hSSchritt\hss}%
    \else\ifnum#1=6
      \advance\ycoord by -1
      \raise\ycoord \Einheit\hbox to0pt{\hskip\xcoord \Einheit
        \hbox{\hskip-2pt \vSSchritt}\hss}%
    \else\ifnum#1=7
      \advance\xcoord by -1
      \advance\ycoord by -1
      \raise\ycoord \Einheit\hbox to0pt{\hskip\xcoord \Einheit
         \DSSchritt\hss}
    \else\ifnum#1=8
      \advance\xcoord by -1
      \advance\ycoord by 1
      \raise\ycoord \Einheit\hbox to0pt{\hskip\xcoord \Einheit
         \dSSchritt\hss}
    \fi\fi\fi\fi
    \fi\fi\fi\fi
  \fi\next}
\def\Koordinatenachsen(#1,#2){\unskip
 \hbox to0pt{\hskip-.5pt\vrule height#2 \Einheit width.5pt depth1 \Einheit}%
 \hbox to0pt{\hskip-1 \Einheit \xcoord#1 \advance\xcoord by1
    \vrule height0.25pt width\xcoord \Einheit depth0.25pt\hss}}
\def\Koordinatenachsen(#1,#2)(#3,#4){\unskip
 \hbox to0pt{\hskip-.5pt \ycoord-#4 \advance\ycoord by1
    \vrule height#2 \Einheit width.5pt depth\ycoord \Einheit}%
 \hbox to0pt{\hskip-1 \Einheit \hskip#3\Einheit 
    \xcoord#1 \advance\xcoord by1 \advance\xcoord by-#3 
    \vrule height0.25pt width\xcoord \Einheit depth0.25pt\hss}}
\def\Gitter(#1,#2){\unskip \xcoord0 \ycoord0 \leavevmode
  \LOOP\ifnum\ycoord<#2
    \loop\ifnum\xcoord<#1
      \raise\ycoord \Einheit\hbox to0pt{\hskip\xcoord \Einheit\Punkt\hss}%
      \advance\xcoord by1
    \repeat
    \xcoord0
    \advance\ycoord by1
  \REPEAT}
\def\Gitter(#1,#2)(#3,#4){\unskip \xcoord#3 \ycoord#4 \leavevmode
  \LOOP\ifnum\ycoord<#2
    \loop\ifnum\xcoord<#1
      \raise\ycoord \Einheit\hbox to0pt{\hskip\xcoord \Einheit\Punkt\hss}%
      \advance\xcoord by1
    \repeat
    \xcoord#3
    \advance\ycoord by1
  \REPEAT}
\def\Label#1#2(#3,#4){\unskip \xdim#3 \Einheit \ydim#4 \Einheit
  \def\lo{\advance\xdim by-.5 \Einheit \advance\ydim by.5 \Einheit}%
  \def\llo{\advance\xdim by-.25cm \advance\ydim by.5 \Einheit}%
  \def\loo{\advance\xdim by-.5 \Einheit \advance\ydim by.25cm}%
  \def\o{\advance\ydim by.25cm}%
  \def\ro{\advance\xdim by.5 \Einheit \advance\ydim by.5 \Einheit}%
  \def\rro{\advance\xdim by.25cm \advance\ydim by.5 \Einheit}%
  \def\roo{\advance\xdim by.5 \Einheit \advance\ydim by.25cm}%
  \def\l{\advance\xdim by-.30cm}%
  \def\r{\advance\xdim by.30cm}%
  \def\lu{\advance\xdim by-.5 \Einheit \advance\ydim by-.6 \Einheit}%
  \def\llu{\advance\xdim by-.25cm \advance\ydim by-.6 \Einheit}%
  \def\luu{\advance\xdim by-.5 \Einheit \advance\ydim by-.30cm}%
  \def\u{\advance\ydim by-.30cm}%
  \def\ru{\advance\xdim by.5 \Einheit \advance\ydim by-.6 \Einheit}%
  \def\rru{\advance\xdim by.25cm \advance\ydim by-.6 \Einheit}%
  \def\ruu{\advance\xdim by.5 \Einheit \advance\ydim by-.30cm}%
  #1\raise\ydim\hbox to0pt{\hskip\xdim
     \vbox to0pt{\vss\hbox to0pt{\hss$#2$\hss}\vss}\hss}%
}
\catcode`\@=13

\catcode`\@=11
\font\tenln    = line10
\font\tenlnw   = linew10

\newskip\Einheit \Einheit=0.5cm
\newcount\xcoord \newcount\ycoord
\newdimen\xdim \newdimen\ydim \newdimen\PfadD@cke \newdimen\Pfadd@cke

\newcount\@tempcnta
\newcount\@tempcntb

\newdimen\@tempdima
\newdimen\@tempdimb

\newdimen\@wholewidth
\newdimen\@halfwidth

\newcount\@xarg
\newcount\@yarg
\newcount\@yyarg
\newbox\@linechar
\newbox\@tempboxa
\newdimen\@linelen
\newdimen\@clnwd
\newdimen\@clnht

\newif\if@negarg

\def\@whilenoop#1{}
\def\@whiledim#1\do #2{\ifdim #1\relax#2\@iwhiledim{#1\relax#2}\fi}
\def\@iwhiledim#1{\ifdim #1\let\@nextwhile=\@iwhiledim
        \else\let\@nextwhile=\@whilenoop\fi\@nextwhile{#1}}

\def\@whileswnoop#1\fi{}
\def\@whilesw#1\fi#2{#1#2\@iwhilesw{#1#2}\fi\fi}
\def\@iwhilesw#1\fi{#1\let\@nextwhile=\@iwhilesw
         \else\let\@nextwhile=\@whileswnoop\fi\@nextwhile{#1}\fi}

\def\thinlines{\let\@linefnt\tenln \let\@circlefnt\tencirc
  \@wholewidth\fontdimen8\tenln \@halfwidth .5\@wholewidth}
\def\thicklines{\let\@linefnt\tenlnw \let\@circlefnt\tencircw
  \@wholewidth\fontdimen8\tenlnw \@halfwidth .5\@wholewidth}
\thinlines

\PfadD@cke1pt \Pfadd@cke0.5pt
\def\PfadDicke#1{\PfadD@cke#1 \divide\PfadD@cke by2 \Pfadd@cke\PfadD@cke \multiply\PfadD@cke by2}
\long\def\LOOP#1\REPEAT{\def\BODY{#1}\ITERATE}
\def\ITERATE{\BODY \let\next\ITERATE \else\let\next\relax\fi \next}
\let\REPEAT=\fi
\def\Punkt{\hbox{\raise-2pt\hbox to0pt{\hss$\ssize\bullet$\hss}}}
\def\DuennPunkt(#1,#2){\unskip
  \raise#2 \Einheit\hbox to0pt{\hskip#1 \Einheit
          \raise-2.5pt\hbox to0pt{\hss$\bullet$\hss}\hss}}
\def\NormalPunkt(#1,#2){\unskip
  \raise#2 \Einheit\hbox to0pt{\hskip#1 \Einheit
          \raise-3pt\hbox to0pt{\hss\twelvepoint$\bullet$\hss}\hss}}
\def\DickPunkt(#1,#2){\unskip
  \raise#2 \Einheit\hbox to0pt{\hskip#1 \Einheit
          \raise-4pt\hbox to0pt{\hss\fourteenpoint$\bullet$\hss}\hss}}
\def\Kreis(#1,#2){\unskip
  \raise#2 \Einheit\hbox to0pt{\hskip#1 \Einheit
          \raise-4pt\hbox to0pt{\hss\fourteenpoint$\circ$\hss}\hss}}

\def\Line@(#1,#2)#3{\@xarg #1\relax \@yarg #2\relax
\@linelen=#3\Einheit
\ifnum\@xarg =0 \@vline
  \else \ifnum\@yarg =0 \@hline \else \@sline\fi
\fi}

\def\@sline{\ifnum\@xarg< 0 \@negargtrue \@xarg -\@xarg \@yyarg -\@yarg
  \else \@negargfalse \@yyarg \@yarg \fi
\ifnum \@yyarg >0 \@tempcnta\@yyarg \else \@tempcnta -\@yyarg \fi
\ifnum\@tempcnta>6 \@badlinearg\@tempcnta0 \fi
\ifnum\@xarg>6 \@badlinearg\@xarg 1 \fi
\setbox\@linechar\hbox{\@linefnt\@getlinechar(\@xarg,\@yyarg)}%
\ifnum \@yarg >0 \let\@upordown\raise \@clnht\z@
   \else\let\@upordown\lower \@clnht \ht\@linechar\fi
\@clnwd=\wd\@linechar
\if@negarg \hskip -\wd\@linechar \def\@tempa{\hskip -2\wd\@linechar}\else
     \let\@tempa\relax \fi
\@whiledim \@clnwd <\@linelen \do
  {\@upordown\@clnht\copy\@linechar
   \@tempa
   \advance\@clnht \ht\@linechar
   \advance\@clnwd \wd\@linechar}%
\advance\@clnht -\ht\@linechar
\advance\@clnwd -\wd\@linechar
\@tempdima\@linelen\advance\@tempdima -\@clnwd
\@tempdimb\@tempdima\advance\@tempdimb -\wd\@linechar
\if@negarg \hskip -\@tempdimb \else \hskip \@tempdimb \fi
\multiply\@tempdima \@m
\@tempcnta \@tempdima \@tempdima \wd\@linechar \divide\@tempcnta \@tempdima
\@tempdima \ht\@linechar \multiply\@tempdima \@tempcnta
\divide\@tempdima \@m
\advance\@clnht \@tempdima
\ifdim \@linelen <\wd\@linechar
   \hskip \wd\@linechar
  \else\@upordown\@clnht\copy\@linechar\fi}

\def\@hline{\ifnum \@xarg <0 \hskip -\@linelen \fi
\vrule height\Pfadd@cke width \@linelen depth\Pfadd@cke
\ifnum \@xarg <0 \hskip -\@linelen \fi}

\def\@getlinechar(#1,#2){\@tempcnta#1\relax\multiply\@tempcnta 8
\advance\@tempcnta -9 \ifnum #2>0 \advance\@tempcnta #2\relax\else
\advance\@tempcnta -#2\relax\advance\@tempcnta 64 \fi
\char\@tempcnta}

\def\Vektor(#1,#2)#3(#4,#5){\unskip\leavevmode
  \xcoord#4\relax \ycoord#5\relax
      \raise\ycoord \Einheit\hbox to0pt{\hskip\xcoord \Einheit
         \Vector@(#1,#2){#3}\hss}}

\def\Vector@(#1,#2)#3{\@xarg #1\relax \@yarg #2\relax
\@tempcnta \ifnum\@xarg<0 -\@xarg\else\@xarg\fi
\ifnum\@tempcnta<5\relax
\@linelen=#3\Einheit
\ifnum\@xarg =0 \@vvector
  \else \ifnum\@yarg =0 \@hvector \else \@svector\fi
\fi
\else\@badlinearg\fi}

\def\@hvector{\@hline\hbox to 0pt{\@linefnt
\ifnum \@xarg <0 \@getlarrow(1,0)\hss\else
    \hss\@getrarrow(1,0)\fi}}

\def\@vvector{\ifnum \@yarg <0 \@downvector \else \@upvector \fi}

\def\@svector{\@sline
\@tempcnta\@yarg \ifnum\@tempcnta <0 \@tempcnta=-\@tempcnta\fi
\ifnum\@tempcnta <5
  \hskip -\wd\@linechar
  \@upordown\@clnht \hbox{\@linefnt  \if@negarg
  \@getlarrow(\@xarg,\@yyarg) \else \@getrarrow(\@xarg,\@yyarg) \fi}%
\else\@badlinearg\fi}

\def\@upline{\hbox to \z@{\hskip -.5\Pfadd@cke \vrule width \Pfadd@cke
   height \@linelen depth \z@\hss}}

\def\@downline{\hbox to \z@{\hskip -.5\Pfadd@cke \vrule width \Pfadd@cke
   height \z@ depth \@linelen \hss}}

\def\@upvector{\@upline\setbox\@tempboxa\hbox{\@linefnt\char'66}\raise
     \@linelen \hbox to\z@{\lower \ht\@tempboxa\box\@tempboxa\hss}}

\def\@downvector{\@downline\lower \@linelen
      \hbox to \z@{\@linefnt\char'77\hss}}

\def\@getlarrow(#1,#2){\ifnum #2 =\z@ \@tempcnta='33\else
\@tempcnta=#1\relax\multiply\@tempcnta \sixt@@n \advance\@tempcnta
-9 \@tempcntb=#2\relax\multiply\@tempcntb \tw@
\ifnum \@tempcntb >0 \advance\@tempcnta \@tempcntb\relax
\else\advance\@tempcnta -\@tempcntb\advance\@tempcnta 64
\fi\fi\char\@tempcnta}

\def\@getrarrow(#1,#2){\@tempcntb=#2\relax
\ifnum\@tempcntb < 0 \@tempcntb=-\@tempcntb\relax\fi
\ifcase \@tempcntb\relax \@tempcnta='55 \or
\ifnum #1<3 \@tempcnta=#1\relax\multiply\@tempcnta
24 \advance\@tempcnta -6 \else \ifnum #1=3 \@tempcnta=49
\else\@tempcnta=58 \fi\fi\or
\ifnum #1<3 \@tempcnta=#1\relax\multiply\@tempcnta
24 \advance\@tempcnta -3 \else \@tempcnta=51\fi\or
\@tempcnta=#1\relax\multiply\@tempcnta
\sixt@@n \advance\@tempcnta -\tw@ \else
\@tempcnta=#1\relax\multiply\@tempcnta
\sixt@@n \advance\@tempcnta 7 \fi\ifnum #2<0 \advance\@tempcnta 64 \fi
\char\@tempcnta}

\def\Diagonale(#1,#2)#3{\unskip\leavevmode
  \xcoord#1\relax \ycoord#2\relax
      \raise\ycoord \Einheit\hbox to0pt{\hskip\xcoord \Einheit
         \Line@(1,1){#3}\hss}}
\def\AntiDiagonale(#1,#2)#3{\unskip\leavevmode
  \xcoord#1\relax \ycoord#2\relax 
      \raise\ycoord \Einheit\hbox to0pt{\hskip\xcoord \Einheit
         \Line@(1,-1){#3}\hss}}
\def\Pfad(#1,#2),#3\endPfad{\unskip\leavevmode
  \xcoord#1 \ycoord#2 \thicklines\ZeichnePfad#3\endPfad\thinlines}
\def\ZeichnePfad#1{\ifx#1\endPfad\let\next\relax
  \else\let\next\ZeichnePfad
    \ifnum#1=1
      \raise\ycoord \Einheit\hbox to0pt{\hskip\xcoord \Einheit
         \vrule height\Pfadd@cke width1 \Einheit depth\Pfadd@cke\hss}%
      \advance\xcoord by 1
    \else\ifnum#1=2
      \raise\ycoord \Einheit\hbox to0pt{\hskip\xcoord \Einheit
        \hbox{\hskip-\PfadD@cke\vrule height1 \Einheit width\PfadD@cke depth0pt}\hss}%
      \advance\ycoord by 1
    \else\ifnum#1=3
      \raise\ycoord \Einheit\hbox to0pt{\hskip\xcoord \Einheit
         \Line@(1,1){1}\hss}
      \advance\xcoord by 1
      \advance\ycoord by 1
    \else\ifnum#1=4
      \raise\ycoord \Einheit\hbox to0pt{\hskip\xcoord \Einheit
         \Line@(1,-1){1}\hss}
      \advance\xcoord by 1
      \advance\ycoord by -1
    \else\ifnum#1=5
      \advance\xcoord by -1
      \raise\ycoord \Einheit\hbox to0pt{\hskip\xcoord \Einheit
         \vrule height\Pfadd@cke width1 \Einheit depth\Pfadd@cke\hss}%
    \else\ifnum#1=6
      \advance\ycoord by -1
      \raise\ycoord \Einheit\hbox to0pt{\hskip\xcoord \Einheit
        \hbox{\hskip-\PfadD@cke\vrule height1 \Einheit width\PfadD@cke depth0pt}\hss}%
    \else\ifnum#1=7
      \advance\xcoord by -1
      \advance\ycoord by -1
      \raise\ycoord \Einheit\hbox to0pt{\hskip\xcoord \Einheit
         \Line@(1,1){1}\hss}
    \else\ifnum#1=8
      \advance\xcoord by -1
      \advance\ycoord by +1
      \raise\ycoord \Einheit\hbox to0pt{\hskip\xcoord \Einheit
         \Line@(1,-1){1}\hss}
    \fi\fi\fi\fi
    \fi\fi\fi\fi
  \fi\next}
\def\hSSchritt{\leavevmode\raise-.4pt\hbox to0pt{\hss.\hss}\hskip.2\Einheit
  \raise-.4pt\hbox to0pt{\hss.\hss}\hskip.2\Einheit
  \raise-.4pt\hbox to0pt{\hss.\hss}\hskip.2\Einheit
  \raise-.4pt\hbox to0pt{\hss.\hss}\hskip.2\Einheit
  \raise-.4pt\hbox to0pt{\hss.\hss}\hskip.2\Einheit}
\def\vSSchritt{\vbox{\baselineskip.2\Einheit\lineskiplimit0pt
\hbox{.}\hbox{.}\hbox{.}\hbox{.}\hbox{.}}}
\def\DSSchritt{\leavevmode\raise-.4pt\hbox to0pt{%
  \hbox to0pt{\hss.\hss}\hskip.2\Einheit
  \raise.2\Einheit\hbox to0pt{\hss.\hss}\hskip.2\Einheit
  \raise.4\Einheit\hbox to0pt{\hss.\hss}\hskip.2\Einheit
  \raise.6\Einheit\hbox to0pt{\hss.\hss}\hskip.2\Einheit
  \raise.8\Einheit\hbox to0pt{\hss.\hss}\hss}}
\def\dSSchritt{\leavevmode\raise-.4pt\hbox to0pt{%
  \hbox to0pt{\hss.\hss}\hskip.2\Einheit
  \raise-.2\Einheit\hbox to0pt{\hss.\hss}\hskip.2\Einheit
  \raise-.4\Einheit\hbox to0pt{\hss.\hss}\hskip.2\Einheit
  \raise-.6\Einheit\hbox to0pt{\hss.\hss}\hskip.2\Einheit
  \raise-.8\Einheit\hbox to0pt{\hss.\hss}\hss}}
\def\SPfad(#1,#2),#3\endSPfad{\unskip\leavevmode
  \xcoord#1 \ycoord#2 \ZeichneSPfad#3\endSPfad}
\def\ZeichneSPfad#1{\ifx#1\endSPfad\let\next\relax
  \else\let\next\ZeichneSPfad
    \ifnum#1=1
      \raise\ycoord \Einheit\hbox to0pt{\hskip\xcoord \Einheit
         \hSSchritt\hss}%
      \advance\xcoord by 1
    \else\ifnum#1=2
      \raise\ycoord \Einheit\hbox to0pt{\hskip\xcoord \Einheit
        \hbox{\hskip-2pt \vSSchritt}\hss}%
      \advance\ycoord by 1
    \else\ifnum#1=3
      \raise\ycoord \Einheit\hbox to0pt{\hskip\xcoord \Einheit
         \DSSchritt\hss}
      \advance\xcoord by 1
      \advance\ycoord by 1
    \else\ifnum#1=4
      \raise\ycoord \Einheit\hbox to0pt{\hskip\xcoord \Einheit
         \dSSchritt\hss}
      \advance\xcoord by 1
      \advance\ycoord by -1
    \else\ifnum#1=5
      \advance\xcoord by -1
      \raise\ycoord \Einheit\hbox to0pt{\hskip\xcoord \Einheit
         \hSSchritt\hss}%
    \else\ifnum#1=6
      \advance\ycoord by -1
      \raise\ycoord \Einheit\hbox to0pt{\hskip\xcoord \Einheit
        \hbox{\hskip-2pt \vSSchritt}\hss}%
    \else\ifnum#1=7
      \advance\xcoord by -1
      \advance\ycoord by -1
      \raise\ycoord \Einheit\hbox to0pt{\hskip\xcoord \Einheit
         \DSSchritt\hss}
    \else\ifnum#1=8
      \advance\xcoord by -1
      \advance\ycoord by 1
      \raise\ycoord \Einheit\hbox to0pt{\hskip\xcoord \Einheit
         \dSSchritt\hss}
    \fi\fi\fi\fi
    \fi\fi\fi\fi
  \fi\next}
\def\Koordinatenachsen(#1,#2){\unskip
 \hbox to0pt{\hskip-.5pt\vrule height#2 \Einheit width.5pt depth1 \Einheit}%
 \hbox to0pt{\hskip-1 \Einheit \xcoord#1 \advance\xcoord by1
    \vrule height0.25pt width\xcoord \Einheit depth0.25pt\hss}}
\def\Koordinatenachsen(#1,#2)(#3,#4){\unskip
 \hbox to0pt{\hskip-.5pt \ycoord-#4 \advance\ycoord by1
    \vrule height#2 \Einheit width.5pt depth\ycoord \Einheit}%
 \hbox to0pt{\hskip-1 \Einheit \hskip#3\Einheit 
    \xcoord#1 \advance\xcoord by1 \advance\xcoord by-#3 
    \vrule height0.25pt width\xcoord \Einheit depth0.25pt\hss}}
\def\Gitter(#1,#2){\unskip \xcoord0 \ycoord0 \leavevmode
  \LOOP\ifnum\ycoord<#2
    \loop\ifnum\xcoord<#1
      \raise\ycoord \Einheit\hbox to0pt{\hskip\xcoord \Einheit\Punkt\hss}%
      \advance\xcoord by1
    \repeat
    \xcoord0
    \advance\ycoord by1
  \REPEAT}
\def\Gitter(#1,#2)(#3,#4){\unskip \xcoord#3 \ycoord#4 \leavevmode
  \LOOP\ifnum\ycoord<#2
    \loop\ifnum\xcoord<#1
      \raise\ycoord \Einheit\hbox to0pt{\hskip\xcoord \Einheit\Punkt\hss}%
      \advance\xcoord by1
    \repeat
    \xcoord#3
    \advance\ycoord by1
  \REPEAT}
\def\Label#1#2(#3,#4){\unskip \xdim#3 \Einheit \ydim#4 \Einheit
  \def\lo{\advance\xdim by-.5 \Einheit \advance\ydim by.5 \Einheit}%
  \def\llo{\advance\xdim by-.25cm \advance\ydim by.5 \Einheit}%
  \def\loo{\advance\xdim by-.5 \Einheit \advance\ydim by.25cm}%
  \def\o{\advance\ydim by.25cm}%
  \def\ro{\advance\xdim by.5 \Einheit \advance\ydim by.5 \Einheit}%
  \def\rro{\advance\xdim by.25cm \advance\ydim by.5 \Einheit}%
  \def\roo{\advance\xdim by.5 \Einheit \advance\ydim by.25cm}%
  \def\l{\advance\xdim by-.30cm}%
  \def\r{\advance\xdim by.30cm}%
  \def\lu{\advance\xdim by-.5 \Einheit \advance\ydim by-.6 \Einheit}%
  \def\llu{\advance\xdim by-.25cm \advance\ydim by-.6 \Einheit}%
  \def\luu{\advance\xdim by-.5 \Einheit \advance\ydim by-.30cm}%
  \def\u{\advance\ydim by-.30cm}%
  \def\ru{\advance\xdim by.5 \Einheit \advance\ydim by-.6 \Einheit}%
  \def\rru{\advance\xdim by.25cm \advance\ydim by-.6 \Einheit}%
  \def\ruu{\advance\xdim by.5 \Einheit \advance\ydim by-.30cm}%
  #1\raise\ydim\hbox to0pt{\hskip\xdim
     \vbox to0pt{\vss\hbox to0pt{\hss$#2$\hss}\vss}\hss}%
}
\catcode`\@=13

\TagsOnRight

\def\BailAA{1}
\def\BaFoAA{2}
\def\CoGuAA{3}
\def\DeWiAA{4}
\def\FrRTAA{5}
\def\GeKrAA{6}
\def\KajiAB{7}
\def\KajiAA{8}
\def\KaNoAA{9}
\def\KawaAA{10}
\def\KratBM{11}
\def\KratBS{12}
\def\RainAA{13}
\def\RoseAA{14}
\def\RoseAB{15}
\def\SchlAR{16}
\def\SlatAC{17}
\def\StanAP{18}
\def\StanBI{19}
\def\WarnAG{20}
\def\WarnAH{21}

\def\AA{1.1}
\def\AB{1.2}
\def\ABb{1.3}
\def\ABc{2.1}
\def\ABd{2.2}
\def\AC{2.3}
\def\AD{2.4}
\def\AE{2.5}
\def\BAa{3.1}
\def\BD{3.2}
\def\BAb{3.3}
\def\BA{3.4}
\def\BDb{3.5}
\def\BDd{3.6}
\def\BB{3.7}
\def\BC{3.8}

\def\TA{1}
\def\EH{2}
\def\TB{3}
\def\TC{4}

\catcode`\@=11
\def\iddots{\mathinner{\mkern1mu\raise\p@\hbox{.}\mkern2mu
    \raise4\p@\hbox{.}\mkern2mu\raise7\p@\vbox{\kern7\p@\hbox{.}}\mkern1mu}}
\catcode`\@=13

\def\fl#1{\left\lfloor#1\right\rfloor}
\def\cl#1{\lceil#1\rceil}

\topmatter 
\title The major index generating function of standard Young
tableaux of shapes of the form ``staircase minus rectangle"
\endtitle 
\author C.~Krattenthaler and M. J. Schlosser
\endauthor 
\affil 
Fakult\"at f\"ur Mathematik, Universit\"at Wien,\\
Oskar-Morgenstern-Platz~1, A-1090 Vienna, Austria.\\
WWW: {\tt http://www.mat.univie.ac.at/\~{}kratt}
WWW: \tt http://www.mat.univie.ac.at/\~{}schlosse
\endaffil
\address Fakult\"at f\"ur Mathematik, Universit\"at Wien,
Oskar-Morgenstern-Platz~1, A-1090 Vienna,\linebreak Austria.\newline
WWW: \tt http://www.mat.univie.ac.at/\~{}kratt,
http://www.mat.univie.ac.at/\~{}schlosse
\endaddress

\thanks Research partially supported 
by the Austrian Science Foundation FWF, grant Z130-N13, grant S9607-N13
(National Research Network
``Analytic Combinatorics and Probabilistic Number Theory"), and
grant SFB F50 (Special Research Program
``Algorithmic and Enumerative Combinatorics")
\endthanks
\subjclass Primary 05A15;
 Secondary 05A19 05E05 33C67
\endsubjclass
\keywords standard Young tableaux, basic hypergeometric series
associated to root systems, elliptic hypergeometric series
associated to root systems
\endkeywords
\abstract 
A specialisation of a transformation formula for multi-dimensional
elliptic hypergeometric series is used to provide compact,
non-determinantal formulae
for the generating function with respect to the major
index of standard Young tableaux of skew shapes of the
form ``staircase minus rectangle". 
\endabstract
\endtopmatter
\document

\subhead 1. Introduction \endsubhead
A {\it standard Young tableau} of skew shape $\la/\mu$,
where $\la=(\la_1,\la_2,\mathbreak\dots,\la_n)$ and
$\mu=(\mu_1,\mu_2,\dots,\mu_n)$ are $n$-tuples of non-negative
integers which are in non-increasing order and satisfy $\la_i\ge\mu_i$
for all $i$, is an arrangement of the
numbers $1,2,\dots,\mathbreak |\la-\mu|$
(where the last quantity denotes the sum of the respective
differences of the integers, $\sum _{i=1} ^{n}(\la_i-\mu_i)$)
of the form
$$\matrix 
&&&\pi_{1,\mu_1+1}&\innerhdotsfor3\after\quad &\pi_{1,\la_1}\\
&&\pi_{2,\mu_2+1}\quad \dots&\pi_{2,\mu_1+1}&\innerhdotsfor2\after\quad 
&\pi_{2,\la_2}\\
&\iddots&&\vdots&&\iddots\\
\pi_{n,\mu_n+1}&\innerhdotsfor3\after\quad &\pi_{n,\la_n}
\endmatrix$$
such that numbers along rows and columns are increasing.
The {\it major index} of a standard Young tableau $T$,
denoted by $\maj(T)$, is defined as the sum over all $i$
such that $i+1$ appears in a lower row in $T$ than $i$.
It is well-known (see e.g\. \cite{\StanBI, Prop.~7.19.11 in
combination with Theorem~7.16.1 and (7.10) with $x_i=q^{i-1}$,
$i=1,2,\dots$, and $y_i=0$ for $i\ge2$})
that the generating function $\sum_T q^{\maj(T)}$,
where the sum runs over all standard Young tableaux of
shape $\la/\mu$, equals
$$
\big[|\la-\mu|]_q!\cdot
\det_{1\le i,j\le n}\left(\frac {1}
{[\la_i-i-\mu_j+j]_q!}\right),
\tag\AA
$$
where $[m]_q!:=[m]_q\,[m-1]_q\cdots[1]_q$ with
$[\al]_q=1+q+q^2+\dots+q^{\al-1}=\frac {1-q^\al} {1-q}$.

\midinsert
$$
\PfadDicke{.5pt}
\Pfad(1,1),222\endPfad
\Pfad(2,2),22\endPfad
\Pfad(3,3),222\endPfad
\Pfad(4,4),22\endPfad
\Pfad(5,5),2\endPfad
\Pfad(0,1),1\endPfad
\Pfad(0,2),11\endPfad
\Pfad(0,3),111\endPfad
\Pfad(2,3),1\endPfad
\Pfad(2,4),11\endPfad
\Pfad(3,5),11\endPfad
\PfadDicke{1pt}
\Pfad(0,0),121212121212\endPfad
\Pfad(0,0),222211122111\endPfad
\hbox{\hskip5cm}
\PfadDicke{.5pt}
\Pfad(1,0),22\endPfad
\Pfad(2,0),22\endPfad
\Pfad(3,0),222222\endPfad
\Pfad(4,1),22222\endPfad
\Pfad(5,2),2222\endPfad
\Pfad(6,3),222\endPfad
\Pfad(7,4),22\endPfad
\Pfad(8,5),2\endPfad
\Pfad(0,1),111\endPfad
\Pfad(2,2),11\endPfad
\Pfad(2,3),111\endPfad
\Pfad(2,4),1111\endPfad
\Pfad(2,5),11111\endPfad
\PfadDicke{1pt}
\Pfad(0,0),11121212121212\endPfad
\Pfad(0,0),22112222111111\endPfad
\hbox{\hskip4cm}
$$
\vskip6pt
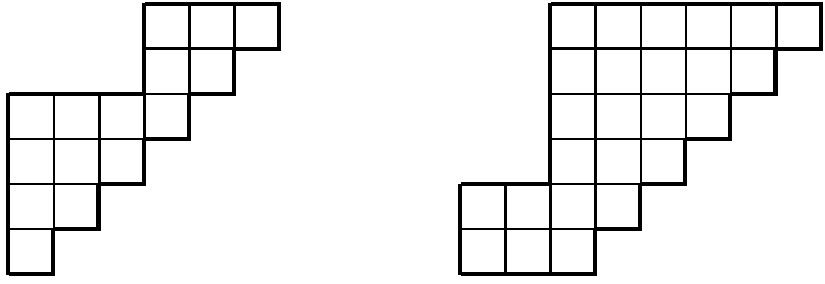
\centerline{\eightpoint Two skew shapes of the form
``staircase minus rectangle"}
\vskip6pt
\centerline{\smc Figure 1}
\vskip0pt
\endinsert

The purpose of this note is to provide formulae
for this major index generating function for 
standard Young tableaux of shapes $\la/\mu$,
where $\la$ is a staircase shape, i.e.,\linebreak
$\la=(N,N-1,\dots,N-n+1)$ for some positive integers $N$ and $n$, 
and $\mu$ is a rectangular shape, i.e.,
$\mu=(m,m,\dots,m,0,\dots,0)$, for some non-negative integer 
$m$ and $r$ repetitions of $m$ (in the sequel we denote
such partitions $\mu$ by $(m^r)$, for short), which are 
(computationally) simpler than the determinantal formula (\AA).
Figure~1 shows the Young diagrams of 
two such shapes according to standard English convention
(cf\. \cite{\StanAP, p.~29}): the diagram on the left represents
the shape $(6,5,4,3,2,1)/(3,3,0,0,0,0)$, and the diagram on the
right represents the shape $(8,7,6,5,4,3)/(2,2,2,2,0,0)$. 
In particular, for $N=n$, the above announced 
formula reduces to a closed form
product formula. To be precise,
we show that the generating function
$\sum_Tq^{\maj(T)}$ for standard Young
tableaux $T$ of shape $(n,n-1,\dots,1)/(m^r)$
(of which the left shape in Figure~1 is an example) equals
$$\multline 
q^{\frac {1} {2}m r(r+m-2n)+\binom {n+1}3}
(1+q)^{\binom n2 -m r }
\left[\tsize\binom {n+1}2-mr\right]_q!\\
\times
\prod_{i=1}^{n}  \frac { [i-1]_{q^2}!} 
        {[2i-1]_q!}
\prod_{i=1}^r  \frac {\big[i-1\big]_{q^2}!\,
   [n+m-r+2i-1]_q!}
   {[m+i-1]_{q^2}!\,[n-m-r+2i-1]_q!},
\endmultline
\tag\AB
$$
a result which was originally (implicitly) obtained
by DeWitt \cite{\DeWiAA, Theorem~V.3} using completely different means.
She proves in fact the stronger result that
a Schur $s$-function of a shape of the form $(n,n-1,\dots,1)/(m^r)$
is the constant multiple of a particular Schur $P$-function.
If this is combined with Kawanaka's product formula for the principal
specialisation of Schur $P$-functions (see \cite{\KawaAA} and
\cite{\RoseAB}), then one obtains the above formula.
Furthermore, for $N=n+1$, we show that the generating function
$\sum_Tq^{\maj(T)}$ for standard Young
tableaux $T$ of shape $(n+1,n,\dots,2)/(m^r)$ equals
$$\multline  
(1+q)^{\binom n2 -(m-1)r}
q^{\frac {1} {2}m r(r+m-2n+2) +r 
\left( 1-n- m  \right) +\binom {n+1}3+   \binom n2 }\\
\times
\left[\tsize\binom {n+2}2-mr-1\right]_q!
\prod_{i=1}^{n} \frac { [i-1]_{q^2}!} 
       { [2i]_q!}
\prod_{i=1}^r  \frac {\big[i-1\big]_{q^2}!\,
    [n+m-r+2i-1]_q!}
   {[m+i-1]_{q^2}!\,[n-m-r+2i]_q!}
\\
\times
\sum_{\ell=0} ^{r}
\frac {(-1)^rq^{2nl_1}} {(1-q^2)^{ r }}
    \bmatrix r\\ {\ell}\endbmatrix_{q^2}
 \frac { \big(q^{-2n};q^2\big)_{\ell}\,
     \big(q^ {n+m-r};q^2\big)_{r-\ell}\,
   \big(q^{n-m-r+1};q^2\big)_{r-\ell}}
    {  \big(q^{n+m-r+1};q^2\big)_{r-\ell}}.
\endmultline\tag\ABb$$
Here, the shifted $q$-factorials are defined by
$(\alpha;q)_k:=(1-\alpha)(1-\alpha q)\cdots(1-\alpha q^{k-1})$ for
$k\ge1$, and $(\alpha ;q)_0:=1$.

In general, if $N=n+s$, where $s$ is a non-negative integer, 
we are able to express the major index generating function
for standard Young
tableaux of shape $(n+s,n+s-1,\dots,s+1)/(m^r)$ as an $\cl{s/2}$-fold 
basic hypergeometric sum, see Theorem~\TA\ in Section~2. 
If $n$ is large compared to $r$ and $s$, then this formula is
computationally superior to the determinantal formula (\AA).

A notable feature of the proof of Theorem~\TA\ that we give here 
is that we require
a basic hypergeometric specialisation of a transformation formula
for multi-dimensional
elliptic hypergeometric series due to Rains and, independently,
Coskun and Gustafson; see Section~3.

\subhead 2. The main result
\endsubhead                
In this section we present our main result, a multi-dimensional
basic hypergeometric series which gives the major index generating
function for standard Young
tableaux of a skew shape that is the difference between a staircase and a
rectangle. 

\proclaim{Theorem \TA}
Let $N,n,m,r$ be non-negative integers with $N\ge n$ and\linebreak 
$N-r+1\ge m$.
If $N-n$ is even, 
the generating function
$\sum_Tq^{\maj(T)}$ for standard Young tableaux $T$ of shape
$(N,N-1,\dots,N-n+1)/(m^r)$ equals
$$\multline   (-1)^{\binom {(N-n)/2}2+\frac {1} {2}r(N-n)}
(1+q)^{\binom n2-\binom {(N-n)/2}2 -m r }
(1-q)^{-\binom {(N-n)/2}2 -r (N-n)}\\
\times
q^{\frac {1} {2}m r(r+m-2n) +\frac {1} {2} r(N-n) 
\left( \frac {1} {2}(N-3n) - m +1 \right) +\binom {n+1}3+ 
(N-n) \left( \binom n2 +\binom {(N-n)/2}2\right)}\\
\times
\frac {\left[\tsize\binom {N+1}2-\binom {N-n+1}2-mr\right]_q!} 
  {\big[r+\frac {N-n-2}2\big]_{q^2}!^{(N-n)/2}\,
  \big[\frac {N+n-2}2\big]_{q^2}!^{(N-n)/2}}
  \frac {\prod_{i=1}^{(N+n)/2} [i-1]_{q^2}!} 
        {\prod_{i=1}^n[N-n+2i-1]_q!}\\
\times
\prod_{i=1}^r  \frac {\big[\frac {N-n} {2}+i-1\big]_{q^2}!\,
   [n+m-r+2i-1]_q!\,\big(q^{n+m-r+2i};q^2\big)_{(N-n)/2}}
   {[m+i-1]_{q^2}!\,[N-m-r+2i-1]_q!}
\\
\times
\sum_{0\le \ell_1<\ell_2<\dots<\ell_{(N-n)/ {2}}\le r+\frac {N-n-2} {2}}
q^{\sum_{i=1}^{(N-n)/2}(N+n-2(2i-1))l_i}
\; \prod_{1\le i<j\le \frac {N-n} {2}}
[\ell_j-\ell_i]_{q^2}^2\kern1.5cm\\
\cdot
\prod_{i=1}^{\frac {N-n}2} \Bigg(
  \bmatrix \frac {N-n-2} {2}+r\\ {\ell_i}\endbmatrix_{q^2}
   \big(q^{2-N-n};q^2\big)_{\ell_i}\,
    \big(q^{n+m-r-2i+1};q^2\big)_{r+i-\ell_i-1}\\
\cdot  \frac {  
      \big(q^{N-m-r-2i+2};q^2\big)_{r+i-\ell_i-1}}
   {\big(q^{N+m-r-2i+2};q^2\big)_{r+i-\ell_i-1}}\Bigg),
\endmultline\tag\ABc$$
while, if $N-n$ is odd, it equals
$$\multline   (-1)^{\binom {(N-n+1)/2}2+\frac {1} {2}r(N-n+1)}
(1+q)^{\binom n2-\binom {(N-n+1)/2}2 -m r }
(1-q)^{-\binom {(N-n+1)/2}2 -r (N-n)}\\
\times
q^{\frac {1} {2}m r(r+m-2n+2) +\frac {1} {2} r(N-n+1) 
\left( \frac {1} {2}(N-3n+1) - m  \right) +\binom {n+1}3+ 
(N-n)  \binom n2 +(N-n+1)\binom {(N-n+1)/2}2}\\
\times
\frac {\left[\tsize\binom {N+1}2-\binom {N-n+1}2-mr\right]_q!} 
  {\big[r+\frac {N-n-1}2\big]_{q^2}!^{(N-n+1)/2}\,
  \big[\frac {N+n-1}2\big]_{q^2}!^{(N-n+1)/2}}
 \frac {\prod_{i=1}^{(N+n+1)/2} [i-1]_{q^2}!} 
       {\prod_{i=1}^n [N-n+2i-1]_q!}\\
\times
\prod_{i=1}^r  \frac {\big[\frac {N-n+1} {2}+i-1\big]_{q^2}!\,
    [n+m-r+2i-1]_q!\,\big(q^{n+m-r+2i+1};q^2\big)_{(N-n-1)/2}}
   {[m+i-1]_{q^2}!\,[N-m-r+2i-1]_q!}
\\
\times
\sum_{0\le \ell_1<\ell_2<\dots<\ell_{(N-n+1)/ {2}}\le r+\frac {N-n-1} {2}}
q^{\sum_{i=1}^{(N-n+1)/2}(N+n+1-2(2i-1))l_i}
\; \prod_{1\le i<j\le \frac {N-n+1} {2}}
[\ell_j-\ell_i]_{q^2}^2\kern0cm\\
\cdot
\prod_{i=1}^{\frac {N-n+1}2} \Bigg(
    \bmatrix \frac {N-n-1} {2}+r\\ {\ell_i}\endbmatrix_{q^2}
 \big(q^{1-N-n};q^2\big)_{\ell_i}\,
     \big(q^ {n+m-r-2i+2};q^2\big)_{r+i-\ell_i-1}\\
\cdot  \frac {
   \big(q^{N-m-r-2i+2};q^2\big)_{r+i-\ell_i-1}}
    {  \big(q^{N+m-r-2i+2};q^2\big)_{r+i-\ell_i-1}}\Bigg).
\endmultline\tag\ABd$$
\endproclaim

\demo{Proof}
According to Formula (\AA), the major index
generating function for standard Young
tableaux which we want to compute is equal to
$$
\left[\tsize\binom {N+1}2-\binom {N-n+1}2-mr\right]_q!
\det_{1\le i,j\le n}\left(
\cases 
\dsize \frac {1} {[N+1-2i-m+j]_q!}&j\le r\\
\dsize \frac {1} {[N+1-2i+j]_q!}&j> r
\endcases
\right).
$$
We now do a Laplace expansion with respect to the first $r$ columns.
In this way we obtain
$$\multline 
\left[\tsize\binom {N+1}2-\binom {N-n+1}2-mr\right]_q!\\
\times
\sum _{1\le k_1<\dots <k_r\le n} ^{}(-1)^{\binom {r+1}2+\sum _{i=1} ^{r}k_i}
\det_{1\le i,j\le r}\left(\frac {1} {[N+1-2k_i-m+j]_q!}\right)\\
\cdot
\underset r+1\le j\le n\to{\det_{1\le i\le n,\,i\notin\{k_1,\dots,k_r\}}}
\left(\frac {1} {[N+1-2i+j]_q!}\right).
\endmultline$$
Using the simple determinant evaluation
$$
\det_{1\le i,j\le s}\left(\frac {1}
{[X_i+j]_q!}\right)=
q^{2\binom {s+1}3+\sum_{i=1}^s(i-1)X_i}
\prod _{i=1} ^{s}\frac {1} {[X_i+s]_q!}
\prod _{1\le i<j\le s} ^{}[X_i-X_j]
$$
(which is readily established by writing $[\al]_q=(1-q^\al)/(1-q)$
for every $q$-integer, factoring out common denominators of rows,
and reduction to a Vandermonde determinant),
the above determinants can be evaluated, so that we arrive at a
multi-dimensional series of basic hypergeometric type:
$$\multline 
\left[\tsize\binom {N+1}2-\binom {N-n+1}2-mr\right]_q!
\kern-6pt\\
\times
\sum _{1\le k_1<\dots <k_r\le n} ^{}(-1)^{\binom {r+1}2+\sum _{i=1} ^{r}k_i}
q^{2\binom {r+1}3+2\binom {n-r+1}3}\kern3cm
\\
\cdot
q^{\sum_{i=1}^r (i-1)(N+1-2k_i-m)+\sum_{i=1}^{n-r} (i-1)(N+1+r)}\\
\cdot
q^{-2\sum_{i=1}^n(i-1)i+2\sum_{i=1}^r (k_i-1)k_i
+2\sum_{i=1}^r i\left(\binom {k_{i+1}}2-\binom {k_i+1}2\right)}
\\
\cdot
\prod _{i=1} ^{n}\frac {1} {[N+n+1-2i]_q!}
\prod _{i=1} ^{r}\frac {[N+n+1-2k_i]_q!} {[N+1-2k_i-m+r]_q!}\\
\cdot
\prod _{1\le i<j\le r} ^{}\big[(N+1-2k_i-m)-(N+1-2k_j-m)\big]_q\\
\cdot
\underset i,j\notin\{k_1,\dots,k_r\}\to
{\prod _{1\le i<j\le n} ^{}}\big[(N+1-2i+r)-(N+1-2j+r)\big]_q\\
=
\left[\tsize\binom {N+1}2-\binom {N-n+1}2-mr\right]_q!\
\prod _{i=1} ^{n}\frac {1} {[N+n+1-2i]_q!}
\prod _{1\le i<j\le n} ^{}[2j-2i]_q
\kern1.5cm
\\
\times
\sum _{1\le k_1<\dots <k_r\le n} ^{}
(-1)^{\binom {r+1}2+\sum _{i=1} ^{r}k_i}
q^{2\binom {r+1}3+2\binom {n-r+1}3
+(N+1-m)\binom r2+(N+1+r)\binom {n-r}2}
\\
\cdot
q^{-4\binom {n+1}3+2r\binom {n+1}2+2\sum_{i=1}^r \binom {k_i}2
-2\sum_{i=1}^r (2i-1)k_i}
\\
\cdot
\prod _{i=1} ^{r}\frac {[N+n+1-2k_i]_q!} {[N+1-2k_i-m+r]_q!}
\prod _{i=1} ^{r}\frac {1} {[2k_i-2]_q!!\,[2n-2k_i]_q!!}
\prod _{1\le i<j\le r} ^{}[2k_j-2k_i]_q^2\\
=
(-1)^{\binom r2}(1+q)^{\binom n2-(n-1)r}(1-q)^{-r(r-1)}
\kern6.5cm\\
\times
q^{2\binom {r+1}3+2\binom {n-r+1}3
+(N+1-m)\binom r2+(N+1+r)\binom {n-r}2
-4\binom {n+1}3+2r\binom {n+1}2-2r^2}
\\
\times
\left[\tsize\binom {N+1}2-\binom {N-n+1}2-mr\right]_q!
\kern5cm\\
\times
\prod _{i=1} ^{n}\frac {[i-1]_{q^2}!} {[N+n+1-2i]_q!}
\prod _{i=1} ^{r}\frac {[N+n-1]_q!} {[n-1]_{q^2}!\,[N-m+r-1]_q!}
\kern.6cm\\
\times
\sum _{0\le k_1<\dots <k_r\le n-1} ^{}
q^{-2\sum_{i=1}^r (2i-1)k_i}
\prod _{1\le i<j\le r} ^{}(1-q^{-2(k_i-k_j)})^2
\kern1cm\\
\cdot
\prod _{i=1} ^{r}\frac {\left(q^{N-m+r-1};q^{-2}\right)_{k_i}\,
\left(q^{N-m+r-2};q^{-2}\right)_{k_i}\,(q^{2n-2};q^{-2})_{k_i}} 
{\left(q^{N+n-1};q^{-2}\right)_{k_i}\,
\left(q^{N+n-2};q^{-2}\right)_{k_i}\,(q^{-2};q^{-2})_{k_i}},
\endmultline\tag\AC$$
where the double $q$-factorials are defined by 
$[2\al]_q!!=[2\al]_q\,[2\al-2]_q\cdots[2]_q$, and, by convention,
$k_{r+1}=n+1$.

At this point, we should note that we may change the range of summation
to 
$$0\le k_1<\dots <k_r\le \fl{\tfrac {1} {2}(N-m+r-1)}.
\tag\AD
$$ 
Indeed, the summand in (\AC) vanishes for $k_r>n-1$ because of the 
term $(q^{2n-2};q^{-2})_{k_i}$, 
unless this incident is ``neutralised" by vanishing
terms in the denominator. However,
since $n\ge r$, we have 
$$N-m+r-1\le N+n-2,
\tag\AE$$ 
except if $m=0$ and $n=r$.
In the latter case, one sees directly that, due to cancellation of
terms, the extension to (\AD) does not change the sum, whereas in the
former case the inequality (\AE) guarantees that the summand still
vanishes for $n-1<k_i\le\frac {1} {2}(N-m+r-1)$.

This last observation makes it possible to apply the transformation formula
in Corollary~\TC\ in Section~3. In order to do so, we have
to distinguish between $N-n$ being even or odd.
If $N-n$ is even, then we have
$$
\frac {(q^{2n-2};q^{-2})_{k_i}} {\left(q^{N+n-2};q^{-2}\right)_{k_i}}
=\frac {\left(q^{N+n-2-2k_i};q^{-2}\right)_{(N-n)/2}} 
{\left(q^{N+n-2};q^{-2}\right)_{(N-n)/2}}.
$$
We then apply (\BDb) with $q$ replaced by $q^{-2}$,
$s=\frac {N-n} {2}$, 
$b=q^{N-m+r-1}$, 
$c=q^{N-m+r-2}$, $d=q^{N+n-2}$,
$f=q^{N+n-1}$ to the sum in (\AC).
After considerable (but routine) manipulation of the arising
expression, the result turns out to equal (\ABc).

If $N-n$ is odd, then we have
$$
\frac {(q^{2n-2};q^{-2})_{k_i}} {\left(q^{N+n-1};q^{-2}\right)_{k_i}}
=\frac {\left(q^{N+n-1-2k_i};q^{-2}\right)_{(N-n+1)/2}} 
{\left(q^{N+n-1};q^{-2}\right)_{(N-n+1)/2}}.
$$
Here, we apply (\BDb) with $q$ replaced by $q^{-2}$, 
$s=\frac {N-n+1} {2}$, 
$b=q^{N-m+r-1}$, $c=q^{N-m+r-2}$, 
$d=q^{N+n-1}$,
$f=q^{N+n-2}$ to the sum in (\AC).
Once again, after lengthy, but routine, manipulation of the arising
expression, we arrive at (\ABd).
\quad \quad \qed
\enddemo

\subhead 3. A transformation formula for elliptic hypergeometric
series and its consequences\endsubhead
In this section we derive the transformation formula 
for multi-dimensional basic hypergeometric series which
is used crucially in the proof of Theorem~\TA\ in the previous section. 
It arises
by specialising and taking appropriate limits in a
transformation formula for multi-dimensional elliptic hypergeometric series.
In order to state this formula, we first need to introduce
``elliptic notation." Given a complex number $p$ with 
$|p|<1$, we define a (rescaled) theta function $\theta(x;p)$
by
$$\theta(x;p)=\prod_{j=0}^\infty(1-p^jx)(1-p^{j+1}/x). $$
Furthermore, out of these ``bricks," we build elliptic
shifted factorials. Namely, fixing another complex parameter, $q$ say,
and a non-negative integer $m$, we set
$$
(a;q,p)_m=\theta(a;p)\,\theta(aq;p)\cdots \theta(aq^{m-1};p),
$$
where the right-hand side is understood as $1$ if $m=0$.
For convenience, we also employ the short notation
$$
(a_1,a_2,\dots,a_k;q,p)_m=
(a_1;q,p)_m\,(a_2;q,p)_m\cdots (a_k;q,p)_m.
$$

The following result 
is a special case of a multi-dimensional ${}_{12}V_{11}$
transformation formula conjectured by Warnaar (let $x=q$ in
\cite{\WarnAG, Conj.~6.1}), which has subsequently been
proven (in more generality) by Rains~\cite{\RainAA, Theorem~4.9} and, 
independently,
by Coskun and Gustafson~\cite{\CoGuAA}. In the present form the
identity has been stated in \cite{\SchlAR, Theorem~3.2}
(with a small typo corrected).

\proclaim{Theorem \EH}
Let $a,b,c,d,e,f$ be indeterminates, let $m$ be a nonnegative integer,
and $r\ge 1$. Then we have
$$\multline
\sum_{0\le k_1<k_2<\dots<k_r\le m}
q^{\sum_{i=1}^r(2i-1)k_i}\prod_{1\le i<j\le r}
\theta(q^{k_i-k_j};p)^2\,\theta(aq^{k_i+k_j};p)^2\\\times
\prod_{i=1}^r\frac{\theta(aq^{2k_i};p)(a,b,c,d,e,f,
\la aq^{2-r+m}/ef,q^{-m};q,p)_{k_i}}
{\theta(a;p)(q,aq/b,aq/c,aq/d,aq/e,aq/f,efq^{r-1-m}/\la,aq^{1+m};q,p)_{k_i}}\\
=\prod_{i=1}^r\frac{(b,c,d,ef/a;q,p)_{i-1}}
{(\la b/a,\la c/a,\la d/a,ef/\la;q,p)_{i-1}}
\kern6cm
\\\times
\prod_{i=1}^r\frac{(aq;q,p)_m\,(aq/ef;q,p)_{m+1-r}\,
(\la q/e,\la q/f;q,p)_{m-i+1}}
{(\la q;q,p)_m\,(\la q/ef;q,p)_{m+1-r}\,
(aq/e,aq/f;q,p)_{m-i+1}}\\\times
\sum_{0\le k_1<k_2<\dots<k_r\le m}
q^{\sum_{i=1}^r(2i-1)k_i}\prod_{1\le i<j\le r}
\theta(q^{k_i-k_j};p)^2\,\theta(\la q^{k_i+k_j};p)^2\\\times
\prod_{i=1}^r\frac{\theta(\la q^{2k_i};p)\,(\la,\la b/a,\la c/a,\la d/a,e,f,
\la aq^{2-r+m}/ef,q^{-m};q,p)_{k_i}}
{\theta(\la;p)\,(q,aq/b,aq/c,aq/d,\la q/e,\la q/f,efq^{r-1-m}/a,
\la q^{1+m};q,p)_{k_i}},
\endmultline\tag\BAa$$
where $\la=a^2q^{2-r}/bcd$.
\endproclaim

We now show that, by suitable specialisation, the above elliptic
hypergeometric transformation formula reduces to the
following transformation formula for multi-dimensional basic hypergeometric
series of different dimensions. 

\proclaim{Corollary \TB}
For all non-negative integers $m$, $r$ and $s$, we have
$$\multline
\sum_{0\le k_1<k_2<\dots<k_r\le m}
q^{\sum_{i=1}^r(2i-1)k_i}
\prod_{1\le i<j\le r} (1-q^{k_i-k_j})^2
\prod_{i=1}^r
\frac{(dq^{k_i};q)_{s}\,(b;q)_{k_i}\,(q^{-m};q)_{k_i}}
{(q;q)_{k_i}\,(f;q)_{k_i}}\\
=
\frac { q^{\binom {r+s}3+\binom {r+1}3+s\binom r2
            - m \binom {r+s}2  }} 
    {f^{\binom r2}\,(q;q)_{r+s-1}^{s-1}}
\prod_{i=1}^r\frac{(b;q)_{i-1}\,(bq^{s+r+i-m-1}/f;q)_{m-r+1}}
{(q^{i-m}/f;q)_{m-i+1}}
\kern1cm\\
\times
\prod_{i=1}^{r+s-1}\frac{(q;q)_{i-1}\,(q;q)_m}
{(q;q)_{m-i}}
\prod_{i=r}^{r+s-1}
\frac {(dq^{1-r}/b;q)_{i}} 
{(q;q)_{r+s-i-1}\,(d;q)_{i-r}\,(fq^{1-r-s}/b;q)_{i}}
\\
\times
\sum_{0\le \ell_1<\ell_2<\dots<\ell_s\le r+s-1}
q^{\sum_{i=1}^s(2i-1)\ell_i}
\prod_{1\le i<j\le s} (1-q^{\ell_i-\ell_j})^2\\
\kern4cm
\cdot
\prod_{i=1}^{s}
\frac {(d;q)_{\ell_i}\,(fq^{1-r-s}/b;q)_{\ell_i}\,(q^{1-r-s};q)_{\ell_i}}
{(q;q)_{\ell_i}\,(dq^{1-r}/b;q)_{\ell_i}\,(q^{-m};q)_{\ell_i}}.\\
\endmultline
\tag\BD$$
\endproclaim

\demo{Proof}
In (\BAa), we let $p=0$, $d\to aq/d$, $f\to aq/f$,
and then $a\to0$. After having performed these substitutions
and limits, we arrive at
$$\multline
\sum_{0\le k_1<k_2<\dots<k_r\le m}
q^{\sum_{i=1}^r(2i-1)k_i}
\prod_{1\le i<j\le r} (1-q^{k_i-k_j})^2\\
\cdot
\prod_{i=1}^r
\frac{(b;q)_{k_i}\,(c;q)_{k_i}\,(e;q)_{k_i}\,(q^{-m};q)_{k_i}}
{(q;q)_{k_i}\,(d;q)_{k_i}\,(f;q)_{k_i}\,(bceq^{2r-m-1}/df)_{k_i}}\\
=\prod_{i=1}^r\frac{(b;q)_{i-1}\,(c;q)_{i-1}\,(eq/f;q)_{i-1}}
{(dq^{1-r}/c;q)_{i-1}\,(dq^{1-r}/b;q)_{i-1}\,
(bceq^r/df;q)_{i-1}}\\\times
\prod_{i=1}^r\frac{(f/e;q)_{m+1-r}\,
(dfq^{1-r}/bc;q)_{m-i+1}}
{(dfq^{1-r}/bce;q)_{m+1-r}\,
(f;q)_{m-i+1}}\\\times
\sum_{0\le k_1<k_2<\dots<k_r\le m}
q^{\sum_{i=1}^r(2i-1)k_i}
\prod_{1\le i<j\le r} (1-q^{k_i-k_j})^2\\
\cdot
\prod_{i=1}^r\frac{
(dq^{1-r}/c;q)_{k_i}\,(dq^{1-r}/b;q)_{k_i}\,(e;q)_{k_i}\,(q^{-m};q)_{k_i}}
{(q;q)_{k_i}\,(d;q)_{k_i}\,(dfq^{1-r}/bc)_{k_i}\,(eq^{r-m}/f)_{k_i}}.
\endmultline
\tag\BAb$$
In this multi-dimensional Whipple transformation we let $e\to\infty$.
This yields the identity
$$\multline
\sum_{0\le k_1<k_2<\dots<k_r\le m}
q^{\sum_{i=1}^r(m-2r+2i)k_i}\left(\frac {df} {bc}\right)^{\sum_{i=1}^rk_i}
\prod_{1\le i<j\le r} (1-q^{k_i-k_j})^2\\
\cdot\prod_{i=1}^r
\frac{(b;q)_{k_i}\,(c;q)_{k_i}\,(q^{-m};q)_{k_i}}
{(q;q)_{k_i}\,(d;q)_{k_i}\,(f;q)_{k_i}}\\
=\left(\frac {dq^{1-r}} {bc}\right)^{\binom r2}
\prod_{i=1}^r\frac{(b;q)_{i-1}\,(c;q)_{i-1}}
{(dq^{1-r}/c;q)_{i-1}\,(dq^{1-r}/b;q)_{i-1}}
\cdot\frac{
(dfq^{1-r}/bc;q)_{m-i+1}}
{(f;q)_{m-i+1}}\\\times
\sum_{0\le k_1<k_2<\dots<k_r\le m}
q^{\sum_{i=1}^r(m-r+2i-1)k_i}f^{\sum_{i=1}^rk_i}
\prod_{1\le i<j\le r} (1-q^{k_i-k_j})^2\\
\prod_{i=1}^r\frac{
(dq^{1-r}/c;q)_{k_i}\,(dq^{1-r}/b;q)_{k_i}\,(q^{-m};q)_{k_i}}
{(q;q)_{k_i}\,(d;q)_{k_i}\,(dfq^{1-r}/bc;q)_{k_i}}.
\endmultline
$$
This identity becomes more elegant if we replace $q$ by $1/q$,
$b$ by $1/b$, $c$ by $1/c$, $d$ by $1/d$, and $f$ by $1/f$, namely
$$\multline
\sum_{0\le k_1<k_2<\dots<k_r\le m}
q^{\sum_{i=1}^r(2i-1)k_i}
\prod_{1\le i<j\le r} (1-q^{k_i-k_j})^2
\prod_{i=1}^r
\frac{(b;q)_{k_i}\,(c;q)_{k_i}\,(q^{-m};q)_{k_i}}
{(q;q)_{k_i}\,(d;q)_{k_i}\,(f;q)_{k_i}}\\
=\left(\frac {bc} {dq^{1-r}}\right)^{r(m-r+1)}
\prod_{i=1}^r\frac{(b;q)_{i-1}\,(c;q)_{i-1}}
{(dq^{1-r}/c;q)_{i-1}\,(dq^{1-r}/b;q)_{i-1}}
\cdot\frac{
(dfq^{1-r}/bc;q)_{m-i+1}}
{(f;q)_{m-i+1}}\\\times
\sum_{0\le k_1<k_2<\dots<k_r\le m}
q^{\sum_{i=1}^r(2i-1)k_i}
\prod_{1\le i<j\le r} (1-q^{k_i-k_j})^2\\
\cdot
\prod_{i=1}^r\frac{
(dq^{1-r}/c;q)_{k_i}\,(dq^{1-r}/b;q)_{k_i}\,(q^{-m};q)_{k_i}}
{(q;q)_{k_i}\,(d;q)_{k_i}\,(dfq^{1-r}/bc;q)_{k_i}}.
\endmultline
\tag\BA
$$
Now we let $c=dq^s$ in this transformation formula, where $s$ is
a non-negative integer. 
Due to the term $(q^{-r-s+1};q)_{k_i}$, on the right-hand side
we are now summing over all $k_1,k_2,\dots,k_r$ with
$0\le k_1<k_2<\dots<k_r\le r+s-1$. Let 
$$\{\ell_1,\ell_2,\dots,\ell_s\}=\{0,1,\dots,r+s-1\}\backslash
\{k_1,k_2,\dots,k_r\},$$
with $\ell_1<\ell_2<\dots<\ell_s$. 
After some simplification, we obtain (\BD).\quad \quad \qed
\enddemo

\remark{Remarks}
(1) The multiple Whipple transformation (\BAb) also follows from letting
$t=q^2$ in \cite{\BaFoAA, Eq.~(5.10)}, followed by 
the principal specialisation formula for the Macdonald polynomial
$P_\lambda(q,t)$.

\medskip
(2) The transformation formula (\BD) is one of the rare examples of
a transformation formula between multi-dimensional 
(basic) hypergeometric series of different dimension.
Other examples that we are aware of are 
\cite{\GeKrAA, Sec.~8}, \cite{\KajiAB}, \cite{\KajiAA},
\cite{\KaNoAA, Thm.~2.2}, \cite{\KratBS, Conjecture in Sec.~1},
\cite{\RoseAA, Thm.~3.1}, and of \cite{\WarnAH, Thm.~4.1}.

\endremark

For our purpose, we shall have need of a more flexible statement
of the transformation formula in Corollary~\TB. More precisely, 
in our application in the proof of Theorem~\TA, we have a pair
$(q^{M};q^{-2})_{k_i}\,(q^{M+1};q^{-2})_{k_i}$ of Pochhammer symbols in the
numerator of our multi-dimensional sum (see (\AC)), and
consequently we would like 
a statement of (\BD) in which the roles of $(b;q)_{k_i}$ and
$(q^{-m};q)_{k_i}$ are interchangeable. Indeed, by a standard
polynomial argument (see e.g\. \cite{\BailAA, Sec.~5.1} or
\cite{\SlatAC, Sec.~2.3.4}), one is able to show the
following ``symmetric" form of Corollary~\TB.

\proclaim{Corollary \TC}
Let one of $b$ and $c$ be of the form $q^{-m}$, where $m$ is a
non-negative integer. Then we have
$$\multline
\sum_{0\le k_1<k_2<\dots<k_r}
q^{\sum_{i=1}^r(2i-1)k_i}
\prod_{1\le i<j\le r} (1-q^{k_i-k_j})^2
\prod_{i=1}^r
\frac{(dq^{k_i};q)_{s}\,(b;q)_{k_i}\,(c;q)_{k_i}}
{(q;q)_{k_i}\,(f;q)_{k_i}}\\
=
(-1)^{\binom {r+s}2}
\frac { q^{ \binom {r+1}3+s\binom r2}} 
    {f^{\binom r2}\,(q;q)_{r+s-1}^{s-1}}
\prod_{i=1}^r\frac{(b;q)_{i-1}\,(q/f;q)_\infty\,(bcq^{s+r+i-1}/f;q)_\infty}
{(cq^i/f;q)_\infty\,(bq^{s+i}/f;q)_\infty}
\kern2cm\\
\times
\prod_{i=1}^{r+s-1}{(q;q)_{i-1}\,(c;q)_i}
\prod_{i=r}^{r+s-1}
\frac {(dq^{1-r}/b;q)_{i}} {(q;q)_{r+s-i-1}\,(d;q)_{i-r}\,(fq^{1-r-s}/b;q)_{i}}
\\
\times
\sum_{0\le \ell_1<\ell_2<\dots<\ell_s\le r+s-1}
q^{\sum_{i=1}^s(2i-1)\ell_i}
\prod_{1\le i<j\le s} (1-q^{\ell_i-\ell_j})^2\\
\cdot
\prod_{i=1}^{s}
\frac {(d;q)_{\ell_i}\,(fq^{1-r-s}/b;q)_{\ell_i}\,(q^{1-r-s};q)_{\ell_i}}
{(q;q)_{\ell_i}\,(dq^{1-r}/b;q)_{\ell_i}\,({c};q)_{\ell_i}}.
\endmultline
\tag\BDb$$
\endproclaim
\demo{Sketch of Proof} 
The assertion is obvious if $c=q^{-m}$, where $m$ is a non-negative integer,
because in that case the identity in (\BDb) is directly equivalent
to (\BD).

On the other hand, let us now suppose that $b$ is of the form
$b=q^{-\beta}$, for some non-negative
integer $\beta$. On the right-hand side of (\BD), we write
$$
\frac {(q;q)_m} {(q;q)_{m-i}}=(q^{m-i+1};q)_i
=(-1)^iq^{im-\binom i2}(q^{-m};q)_i
$$
and
$$
\frac{(bq^{s+r+i-m-1}/f;q)_{m-r+1}}
{(q^{i-m}/f;q)_{m-i+1}}=
\frac
{(bq^{s+r+i-m-1}/f;q)_\infty\,(q/f;q)_\infty}
{(bq^{s+i}/f;q)_\infty\,(q^{i-m}/f;q)_\infty}.
$$
Since we are currently assuming that $b=q^{-\beta}$ for some 
non-negative integer $\beta$, the above relation can be 
rewritten in the form
$$
\frac{(bq^{s+r+i-m-1}/f;q)_{m-r+1}}
{(q^{i-m}/f;q)_{m-i+1}}=
\frac{(q/f)_{s+i-\beta-1}} 
{(q^{i-m}/f;q)_{s+r-\beta-1}}
.
\tag\BDd
$$
Consequently,
the expression in (\BDd) is rational in $c=q^{-m}$, and therefore as well
the expressions on both sides of (\BDb). Under a fixed choice
of the parameters $r,s,b=q^{-\beta}$, one also sees that the degrees in $c$ of
numerator and denominator of these rational functions are bounded.
Now comes the ``polynomial argument": Identity~(\BDb) holds for
infinitely many $c$'s (namely for all $c$'s of the form $c=q^{-m}$,
where $m$ is a non-negative integer),
hence it must hold for arbitrary $c$, that is, for $c$ 
considered as an indeterminate. This finishes the proof of this
corollary.\quad \quad \qed
\enddemo

For the convenience of the reader, we state the special cases
of (\BDb) where $s=0$ and $s=1$ explicitly below.
Namely, if we let $s=0$, then (\BDb) reduces to
$$\multline
\sum_{0\le k_1<k_2<\dots<k_r\le m}
q^{\sum_{i=1}^r(2i-1)k_i}
\prod_{1\le i<j\le r} (1-q^{k_i-k_j})^2
\prod_{i=1}^r
\frac{(b;q)_{k_i}\,(q^{-m};q)_{k_i}}
{(q;q)_{k_i}\,(f;q)_{k_i}}\\
=(-1)^{\binom r2}q^{\binom r2(m-r+1)}b^{r(m-r+1)}
\prod_{i=1}^r\frac{(b;q)_{i-1}\,(q;q)_m\,(fq^{1-i}/b;q)_{m-r+1}\,(q;q)_{i-1}}
{(f;q)_{m-i+1}\,(q;q)_{m-i+1}},
\endmultline
\tag\BB$$
where $m$ is a non-negative integer,
while for $s=1$ Identity~(\BDb) reduces to
$$\multline
\sum_{0\le k_1<k_2<\dots<k_r}
q^{\sum_{i=1}^r(2i-1)k_i}
\prod_{1\le i<j\le r} (1-q^{k_i-k_j})^2
\prod_{i=1}^r
\frac{(1-dq^{k_i})\,(b;q)_{k_i}\,(c;q)_{k_i}}
{(q;q)_{k_i}\,(f;q)_{k_i}}\\
= 
(-1)^{\binom {r+1}2}
\frac { q^{ \binom {r+1}3+\binom r2}} 
    {f^{\binom r2}}
\prod_{i=1}^r\frac{(q;q)_{i-1}\,(b;q)_{i-1}\,(c;q)_i\,
               (q/f;q)_\infty\,(bcq^{r+i}/f;q)_\infty}
{(cq^i/f;q)_\infty\,(bq^{i+1}/f;q)_\infty}
\kern2cm\\
\times
\sum_{\ell=0}^r
q^{\ell}
\frac {(d;q)_{\ell}\,(dq^{1-r+\ell}/b;q)_{r-\ell}\,(q^{-r};q)_{\ell}}
{(q;q)_{\ell}\,(fq^{-r+\ell}/b;q)_{r-\ell}\,({c};q)_{\ell}},
\endmultline
\tag\BC$$
where one of $b$ or $c$ is of the form $q^{-m}$, where $m$ is a 
non-negative integer.
We point out that (\BB) has been explicitly stated earlier in
\cite{\KratBM, Theorem~6}, where two different proofs were
provided: one proceeded by showing that this formula
comes from the principal specialisation of the obvious expansion
of a rectangular Schur function in two sets of variables,
while the other derived it by specialising a multi-dimensional
$q$-beta integral formula due to Evans and Kadell.


\Refs

\ref\no \BailAA\by W. N. Bailey \yr 1935 \book Generalized
hypergeometric series\publ Cambridge University Press\publaddr
Cambridge\endref 

\ref\no \BaFoAA\by T.~H.~Baker and P.~J.~Forrester\yr 1999 \paper
Transformation formulas for multivariable basic hypergeometric series\jour
Methods Appl\. Anal\.\vol 6\pages 147--164\endref 

\ref\no \CoGuAA\by H.~Coskun and R.~A.~Gustafson\paper
Well-poised Macdonald functions $W_\la$ and Jackson coefficients
$\omega_\la$ on $BC_n$\inbook Proceedings of the workshop on
Jack, Hall--Littlewood and Macdonald polynomials\eds
V.~B.~Kuznetsov and S.~Sahi\publ
Contemp\. Math., vol.~417\publaddr Amer\. Math\. Soc., Providence, RI
\yr 2006\pages 127--155\endref

\ref\no \DeWiAA\by E. DeWitt\book Identities relating Schur $s$-functions 
and $Q$-functions \publ Ph.D. thesis\publaddr University of Michigan,
Ann Arbor\yr 2012\finalinfo available at
{\tt http://deepblue.lib.umich.edu/handle/2027.42/93841}\endref

\ref\no \FrRTAA\by J. S. Frame, G. B. Robinson and R. M. Thrall \yr
1954 \paper The hook graphs of the symmetric group\jour
Canad\. J. Math\.\vol 6\pages 316--325\endref 

\ref\no \GeKrAA\by I. M. Gessel and C. Krattenthaler \yr 1997 \paper
Cylindric partitions\jour Trans\. Amer\. Math\. Soc\.\vol 349\pages
429--479\endref 

\ref\no \KajiAB\by Y. Kajihara \yr 2001 \paper Some remarks on
multiple Sears transformations\jour Contemp\. Math\.\vol 291\pages
139--145\endref

\ref\no \KajiAA\by Y. Kajihara \yr 2004 \paper Euler transform
formula for multiple basic hypergeometric series of type $A$ and some
applications\jour Adv\. Math\.\vol 187\pages 53--97\endref 

\ref\no \KaNoAA\by Y.~Kajihara and M.~Noumi \yr 2003 \paper
Multiple elliptic hypergeometric series. An approach from the Cauchy
determinant\jour Indag\. Math\. (N.S.) \vol 14\pages 395--421\endref 

\ref\no \KawaAA\by  N. Kawanaka\paper A $q$-series identity involving 
Schur functions and related topics\jour Osaka J. Math\.\vol 36\yr 1999
\pages 157--176\endref

\ref\no \KratBM\by C. Krattenthaler \yr 2000 \paper Schur function 
identities and the number of perfect matchings of holey Aztec 
rectangles\jour Contemporary Math\. \vol 254\pages 335--350\endref

\ref\no \KratBS\by C. Krattenthaler \yr 2001 \paper Proof of a
summation formula for an $\tilde A_n$ basic hypergeometric series
conjectured by Warnaar\inbook ``$q$-Series with Applications to
Combinatorics, Number Theory, and Physics," Urbana--\linebreak Champaign,
Oct.~26--28, 2000\eds B.~C.~Berndt, K.~Ono\publ Contemporary Math.,
vol.~291, Amer\. Math\. Soc\.\publaddr Providence, R.I.\pages
153--161\endref 

\ref\no \RainAA\by E.~Rains\paper
$BC_n$-symmetric abelian functions\jour Duke Math\. J.\vol
135\yr 2006\pages 99--180\endref

\ref\no \RoseAA\by H.~Rosengren\yr 2006\paper New transformations for
elliptic hypergeometric series on the root system $A_n$\jour
Ramanujan J\. \vol 12\pages 155--166\endref 

\ref\no \RoseAB\by H.~Rosengren\yr \paper Schur $Q$-polynomials, 
multiple hypergeometric series and enumeration of marked shifted tableaux
\jour J. Combin. Theory Ser.~A \vol 115 \yr 2008\pages 376--406\endref 

\ref\no \SchlAR\by M. J. Schlosser \paper Elliptic enumeration of 
nonintersecting lattice paths
\jour J. Combin\. Theory Ser.~A \vol 114\yr 2007\pages 505--521\endref

\ref\no \SlatAC\by L. J. Slater \yr 1966 \book Generalized
hypergeometric functions\publ Cambridge University Press\publaddr
Cambridge\endref 

\ref\no \StanAP\by R. P. Stanley \yr 1986 \book Enumerative
Combinatorics\bookinfo vol.~1\publ Wadsworth \& Brooks/Cole\publaddr
Pacific Grove, California\finalinfo reprinted by Cambridge University
Press, Cambridge, 1998\endref 

\ref\no \StanBI\by R. P. Stanley \yr 1999 \book Enumerative 
Combinatorics\bookinfo vol.~2\publ Cambridge University 
Press\publaddr Cambridge\endref

\ref\no \WarnAG\by S.~O.~Warnaar\paper
Summation and transformation formulas for elliptic hypergeometric series
\jour Constr\. Approx\.\vol 18\yr 2002\pages 479--502\endref

\ref\no \WarnAH\by S.~O.~Warnaar\paper The $sl_3$ Selberg
integral\jour Adv\. Math\.\vol 224\yr 2010\pages 499--524\endref

\endRefs
\enddocument